\documentclass{article}

\usepackage{arxiv}

\usepackage[utf8]{inputenc} % allow utf-8 input
\usepackage[T1]{fontenc}    % use 8-bit T1 fonts
\usepackage{hyperref}       % hyperlinks
\usepackage{url}            % simple URL typesetting
\usepackage{booktabs}       % professional-quality tables
\usepackage{amsfonts}       % blackboard math symbols
\usepackage{nicefrac}       % compact symbols for 1/2, etc.
\usepackage{microtype}      % microtypography
\usepackage{lipsum}		% Can be removed after putting your text content
\usepackage{graphicx}
\usepackage{natbib}
\usepackage{doi}
\usepackage{amsmath,amssymb,amsfonts,amsthm}
\usepackage{graphicx}
\usepackage{subcaption}
\usepackage{fancyhdr}
\usepackage{algorithm}
\usepackage{algpseudocode}
\usepackage{comment}
\usepackage{xcolor}

\usepackage{mathtools}
\usepackage{enumitem}
\usepackage{multirow}   % per \multirow nelle tabelle

\newtheorem{theorem}{Theorem}[section]
\newtheorem{lemma}[theorem]{Lemma}
\newtheorem{proposition}[theorem]{Proposition}
\newtheorem{corollary}[theorem]{Corollary}

\newtheorem{remark}[theorem]{Remark}

\title{Limited-Memory LRSGA: An Iterative Method for Computing Nash Equilibria in Competitive Optimization Problems}

\author{Katherine Rossella Foglia%
  \thanks{Department of Mathematics and Computer Science, University of Calabria, Ponte P. Bucci, 30B, Arcavacata di Rende (CS), 87036, Italy. \texttt{katherine.foglia@unical.it}} \\
	Department of Mathematics and Computer Science\\
	University of Calabria\\
	Rende (CS), 87036, Italy \\
	\texttt{katherine.foglia@unical.it} \\
	\And
	Francesco Sergio Pisani%
    \thanks{Institute for High Performance Computing and Networking, Italian National Research Council, Via P. Bucci 8/9C,Arcavacata di Rende (CS), 87036, Italy. \texttt{francescosergio.pisani@icar.cnr.it}} \\
	Institute for High Performance Computing and Networking\\
	Italian National Research Council\\
	Rende (CS), 87036, Italy \\
	\texttt{francescosergio.pisani@icar.cnr.it} \\
	\And
	Vittorio Colao%
    \thanks{Department of Mathematics and Computer Science, University of Calabria, Ponte P. Bucci, 30B, Arcavacata di Rende (CS), 87036, Italy
    \texttt{vittorio.colao@unical.it}} \\
	Department of Mathematics and Computer Science\\
	University of Calabria\\
	Rende (CS), 87036, Italy \\
	\texttt{vittorio.colao@unical.it} \\
}

\date{}

\renewcommand{\shorttitle}{\textit{arXiv} Template}

\hypersetup{
pdftitle={A template for the arxiv style},
pdfsubject={q-bio.NC, q-bio.QM},
pdfauthor={David S.~Hippocampus, Elias D.~Striatum},
pdfkeywords={First keyword, Second keyword, More},
}

\begin{document}
\maketitle

\begin{abstract}
We introduce LM-LRSGA, a limited-memory variant of Low-Rank Symplectic Gradient Adjustment (LRSGA) for differentiable games. 
It is an iterative scheme for approximating Nash equilibria with first-order-like cost while retaining the stabilizing effect of symplectic second-order corrections via low-rank information.
%It is an iterative scheme for approximating Nash equilibria at first-order cost while retaining the stabilizing benefits of second-order information. 
By storing only a limited history of curvature pairs, LM-LRSGA is well suited to high-parameter competitive models such as GANs.
In particular, we provide a per-iteration spectral stability condition for LRSGA near Nash equilibria, a limited‑memory implementation (LM‑LRSGA) based on adapted two‑loop recursions together with a local convergence analysis for fixed history length, and an empirical evaluation on GAN training on MNIST and Fashion-MNIST, including spectral diagnostics of the training dynamics.
\end{abstract}

% keywords can be removed

\section{Introduction}

Optimization in multi-agent and adversarial settings is a central challenge in modern machine learning. Examples include Generative Adversarial Networks (GANs), adversarial training for robustness, and multi-agent reinforcement learning, where two or more agents optimize distinct, often conflicting objectives.

These problems can be viewed as differentiable games: each agent minimizes its own objective. Despite the similarity to classical optimization, the coupled objectives and the lack of a single potential make the dynamics fundamentally different. Even in the two-player case, gradient methods may cycle, oscillate, or diverge rather than converge to equilibrium \cite{arjovsky2017wasserstein, wiatrak2019stabilizing, balduzzi2018mechanics}.
 
Formally, for the two-player setting, we denote by $x\in\mathbb{R}^m$ and $y\in\mathbb{R}^n$ the strategy variables for Players~1 and~2, and by $f,g:\mathbb{R}^{m+n}\to\mathbb{R}$ the corresponding  payoff functions.

A Nash equilibrium $(x^*,y^*)$ is a pair where neither player can locally improve its objective, holding the other fixed. For differentiable $f$ and $g$, the first-order conditions are
\begin{equation}
\partial_x f(x^*,y^*) = 0, 
\qquad 
\partial_y g(x^*,y^*) = 0.
\end{equation}

Stacking $w=(x,y)\in\mathbb{R}^{m+n}$, the game gradient is defined as
\[
F(w)=
\begin{pmatrix}
\partial_x f(w)\\[3pt]
\partial_y g(w)
\end{pmatrix},
\]
so that equilibria satisfy $F(w^*)=0$. This generalizes single-objective optimization to a vector field, but here $F$ is typically non-conservative and non-monotone, so standard convex tools do not apply in general.

Understanding the underlying geometry is therefore crucial. The Jacobian of $F$, i.e., the game Hessian,
\[
H(w) =
\begin{pmatrix}
D(\partial_{x}f(w)) \\
D(\partial_{y}g(w))
\end{pmatrix} =
\begin{pmatrix}
\partial^2_{xx}f(w) & \partial^2_{xy}f(w)\\
\partial^2_{yx}g(w) & \partial^2_{yy}g(w)
\end{pmatrix},
\]
admits the decomposition $H=S+A$ with symmetric $S=\tfrac12(H+H^\top)$ and antisymmetric $A=\tfrac12(H-H^\top)$. The part $S$ drives motion toward equilibria, while $A$ induces rotational or Hamiltonian-like components that prevent convergence. Prior work \cite{balduzzi2018mechanics, gemp2018global} shows that unmitigated rotations cause simultaneous gradient descent to behave as a non-dissipative dynamical system; even in the bilinear game $\Phi(x,y)=x^\top y$,
\[
x_{n+1} = x_n - \eta \,\partial_x f(x_n,y_n), 
\quad 
y_{n+1} = y_n - \eta \,\partial_y g(x_n,y_n)
\]
generates limit cycles (see also \cite{schafer2019competitive}). These rotations arise from antisymmetric interactions between players’ gradients and behave like a skew-symmetric force field. In continuous time, dynamics resemble Hamiltonian systems: energy is conserved, producing a potential (gradient) component toward equilibria and a rotational flow around them. Hence, restoring convergence requires damping or correcting the antisymmetric part.

This insight motivated algorithms with second-order corrections, symplectic damping, or predictive steps. The Symplectic Gradient Adjustment (SGA) \cite{balduzzi2018mechanics} preconditions updates with $(I-\tau A(w))$, counteracting rotation and stabilizing trajectories near equilibria. The Optimistic Gradient method~\cite{daskalakis2018limit} and the Extragradient technique~\cite{korpelevich1976extragradient, gidelvariational2019} exploit the variational structure of the problem to capture and correct the cycling behaviour.
Finally, Competitive Gradient Descent (CGD)~\cite{schafer2019competitive} interprets game optimization through a local bilinear model whose Nash equilibrium defines a stable update.
However, many of these methods are costly: explicit or implicit mixed second-order derivatives scale poorly. Even for moderate networks, forming $\partial^2_{xy}f$ and $\partial^2_{yx}g$ is impractical. Quasi-Newton approximations mitigate this: BFGS/L-BFGS \cite{nocedal1980updating, liu1989limited} maintain low-rank curvature via secant updates, achieving near-second-order behavior at first-order cost. Building on this, Low-Rank SGA (LRSGA) \cite{lrsga} replaces mixed derivatives with low-rank secant updates, recovering much of SGA’s stabilization without full Hessian storage. Yet the low-rank factors still require memory that grows with problem dimension, which can become prohibitive for large models.

This paper extends LRSGA by adding a spectral stability result, an $\varepsilon$-regularized mixed-block formulation, and a limited-memory variant, \emph{Limited-Memory Low-Rank Symplectic Gradient Adjustment} (LM-LRSGA), with local convergence analysis. The method combines symplectic gradient corrections with L-BFGS-style recursions: instead of storing dense factors, it keeps a short history of curvature pairs capturing essential information from recent iterates. Modified two-loop recursion algorithms compute the LRSGA products implicitly, avoiding matrix formation or storage. Thus LM-LRSGA preserves LRSGA’s stabilizing effect while using linear memory and first-order-like per-step complexity, making it suitable for large-scale adversarial/stochastic settings (e.g., GANs with millions of parameters).

Unlike L-BFGS, which aims to approximate the full Hessian structure, LM-LRSGA targets only the antisymmetric component responsible for rotational dynamics. As a result, the method avoids the estimation of the full Hessian information, thereby reducing the amount of second-order approximation required.

We validate the approach on GAN training over MNIST and Fashion-MNIST, comparing LM-LRSGA against the Adam \cite{adam} optimizer. Empirically, LM-LRSGA yields more stable training 
dynamics and faster convergence, with smooth and synchronized loss trajectories for both generator and discriminator. Quantitatively, it achieves lower Fréchet Inception Distance (FID) scores than Adam while maintaining comparable computational cost. Sensitivity analyses show robustness to learning rate, history length, and damping factors, highlighting suitability for large-scale and stochastic regimes.

Beyond computation, our method supports sustainable ML by avoiding explicit second-order computation and reducing memory footprint.
%, LM-LRSGA lowers training time, suggesting efficiency and environmental friendliness

The rest of the paper is organized as follows.

Section~\ref{sec:background} reviews competitive optimization and recalls SGA and LRSGA.Section~3 establishes a spectral stability condition for LRSGA near Nash equilibria. Section~\ref{sec:lmlrsga} presents LM-LRSGA, detailing its two-loop recursion algorithms and update rules.
Section~5 proves local linear convergence of LM-LRSGA under fixed history length and sufficiently small step sizes.Section~\ref{sec:experiments} reports experiments on GANs trained on MNIST and Fashion-MNIST, and presents sensitivity analyses. 
Section~7 complements the empirical results with spectral and stability diagnostics (Jacobian spectrum estimation and Welch-PSD refinement).
Section~\ref{sec:conclusion} concludes and outlines future work, including extending LM-LRSGA to larger-scale stochastic and distributed training, to noncooperative games with more than two players, and quantifying the associated computational and carbon footprint.

\section{Preliminaries and Mathematical Framework}
\label{sec:background}
Let \(X \subset \mathbb{R}^m\) and \(Y \subset \mathbb{R}^n\) be nonempty, convex and compact subsets, and \(\Omega \subset \mathbb{R}^{m+n}\) open and convex such that \(X\times Y \subset \Omega\). Let \(f,g \in C^3(\Omega,\mathbb{R})\) be two objective functions.
A two-player competitive optimization problem, where Player~1 (with strategy $x$) and Player~2 (with strategy $y$) aim to minimize their respective objective functions $f$ and $g$, can be written as:

\begin{equation} \label{problemacomp}
    \min_{x \in X} f(x,y) \qquad \min_{y \in Y} g(x,y).
\end{equation}

Solutions to problem \ref{problemacomp} are called Nash equilibria (NE) and denoted by $(x^*,y^*)$.
For a shorter notation, we write $w = (x, y) \in X \times Y$.

The first-order stationarity and the second-order necessary conditions for a local NE $w^* = (x^*, y^*)$ are
\begin{equation}\label{equiHess1NE}
F(w^*)=0, \qquad \partial_{xx}^2 f(w^*)\succeq 0, \qquad \partial_{yy}^2 g(w^*)\succeq 0.
\end{equation}
That is, the game gradient vanishes at $w^*$ and each player’s second order derivative with respect to its own variable is positive semidefinite.

%Moreover, to guarantee that the Nash equilibrium $w^*$ is isolated and locally stable, we assume that the game Jacobian $H(w^*)$ is nonsingular and that its symmetric part $S(w^*)$ is uniformly positive definite, i.e., there exists $\mu>0$ such that $S(w^*)\succeq \mu I$.

We remark that the game Hessian matrix admits the symmetric/antisymmetric decomposition $H(x,y)=S(x,y)+A(x,y)$:
\[
S(x,y) = \frac{H(x,y) + H(x,y)^{\top}}{2} \quad A(x,y) = \frac{H(x,y)-H(x,y)^{\top}}{2},
\]
and, in order to guarantee local stability in our analysis, we adopt the sufficient assumption that the symmetric part satisfies $S(w^*)\succeq \mu I$ for some $\mu>0$.
This condition provides a uniform coercive (dissipative) component in the linearized game dynamics, which underpins
discrete-time contraction bounds for the preconditioned iteration.
In practice, such coercivity is consistent with (and may be reinforced by) common quadratic regularization terms.

As reported in the literature \cite{balduzzi2018mechanics}, the antisymmetric component \(A(x,y)\) is directly responsible for rotational/oscillatory dynamics during training, which can delay or even prevent rapid convergence to equilibria.

To counteract such phenomena, the \emph{Symplectic Gradient Adjustment (SGA)} method introduced in \cite{balduzzi2018mechanics} modifies the descent direction by subtracting the rotational component induced by \(A\). Its update rule is

\[
w_{k+1} = w_k - \eta\,(I - \tau\,A(w_{k}))\,F(w_k),
\]
where $I$ represents the identity matrix in the appropriate Euclidean space, \(w_k:=(x_k,y_k)\) the iterate $k$, \(\eta>0\) the stepsize, \(\tau>0\) controls the strength of the antisymmetric correction and \(A(w_{k})\) has the block form:

\begin{equation}
    \label{SGA}
A(w_{k})= 
\begin{pmatrix} 
0_m &  \frac{1}{2}\, ( \partial_{xy}^2 f(w_k) - \partial_{yx}^2 g(w_k)^{\top})\\ 
 \frac{1}{2}\, (\partial_{yx}^2 g(w_k) - \partial_{xy}^2 f(w_k)^{\top}) & 0_n 
\end{pmatrix}.
\end{equation}

As outlined in the introduction, our recent variant—\emph{Low-Rank Symplectic Gradient Adjustment (LRSGA)}—proposed in \cite{lrsga}, addresses the high computational cost of explicitly forming the mixed second derivatives \(\partial^2_{xy} f(x,y)\) and \(\partial^2_{yx} g(x,y)\) by employing Broyden’s least–change secant updates \cite{broyden1965class}.
In particular, the employed updates are the following, 
\begin{align}
    \label{eq:mu-nu-interif}
    \mu_{k+1} &= \mu_k +
    \frac{(\partial_x f(w_{k+1}) - \partial_x f(w_k) - \mu_k (w_{k+1} - w_k))(w_{k+1} - w_k)^T}
    {(w_{k+1} - w_k)^T (w_{k+1} - w_k)}, \\
    \label{eq:mu-nu-interig}
    \nu_{k+1} &= \nu_k
    + \frac{(\partial_y g(w_{k+1}) - \partial_y g(w_k) - \nu_k (w_{k+1} - w_k))(w_{k+1} - w_k)^T}
    {(w_{k+1} - w_k)^T (w_{k+1} - w_k)},
\end{align}
obtaining \(\mu_k \approx D(\partial_{x}f(w))=\bigl(\partial_{xx}^2 f \,\big|\, \partial_{xy}^2 f\bigr)\in \mathbb{R}^{m\times(m+n)}\) and \(\nu_k \approx D(\partial_{y}g(w))=\bigl(\partial_{yx}^2 g \,\big|\, \partial_{yy}^2 g\bigr)\in \mathbb{R}^{n\times(m+n)}\).
Henceforth we retain only the mixed–derivative block–columns, denoted \([\mu_k]_2\) and \([\nu_k]_1\), which approximate \(\partial_{xy}^2 f\) and \(\partial_{yx}^2 g\), respectively.
\\
The LRSGA iterative formula is then the following:

\begin{align}
\label{lrsga-alpha}
w_{k+1} = w_k - \eta\,(I - \tau\,\alpha_k)\,F(w_k),
\end{align}
where
\begin{align}
\label{alphak}
\alpha_k =
\begin{pmatrix}
0 & \frac{1}{2}([\mu_k]_2 - [\nu_k]_1^{\top}) \\[6pt]
\frac{1}{2}([\nu_k]_1 - [\mu_k]_2^{\top}) & 0
\end{pmatrix}\approx A(w_k).
\end{align}

This work introduces a \emph{Limited-Memory LRSGA (LM–LRSGA)} variant that preserves the theoretical framework of LRSGA while avoiding the explicit storage of the matrices \(\mu_k\) and \(\nu_k\). The implementation follows a \emph{two-loop recursion} scheme in the spirit of L-BFGS \cite{nocedal1980updating,liu1989limited}, leading to substantial memory savings and making the method suitable for large-scale competitive settings. Further details and evidence are provided in the subsequent sections.

\section{Spectral Analysis of LRSGA}
We establish a per–iteration spectral condition showing how the skew correction $(I-\tau\alpha_k)$ stabilizes the
linearized dynamics at a Nash equilibrium. The key step is an identity for the symmetric part of $(I-\tau K)H^*$ when
$K$ is skew (Lemma~\ref{lem:sym-generalK}); specializing $K=\alpha_k$ yields a spectral condition for the Jacobian of the k-frozen LRSGA map defined by $T^{\mathrm{LRSGA}}_{k}(w)=w-\eta\,(I-\tau\alpha_k)\,F(w)$ (Proposition~\ref{prop:LRSGA-Jac}). Throughout, $\|\cdot\|_2$ denotes the
spectral norm and $\mathrm{sym}(B)$ denote the symmetric part of a generic matrix B.

\begin{lemma}
\label{lem:sym-generalK}
Let $H^*=DF(w^*)=S^*+A^*$ with $S^*=S^{*\top}\succ0$ and $A^*=-A^{*\top}$.
For every skew matrix $K=-K^\top$ and $\tau>0$,
\begin{equation}
\mathrm{sym}\!\big((I-\tau K)H^*\big)
= S^* + \frac{\tau}{2}\big(S^*K-KS^*\big)\;-\;\frac{\tau}{2}\big(KA^*+A^*K\big).
\label{eq:sym-generalK}
\end{equation}
\end{lemma}

\begin{proof}
Using $\mathrm{sym}(M)=\tfrac12(M+M^\top)$ and the skew-symmetry $K^\top=-K$ (hence $(I-\tau K)^\top=I+\tau K$), we have
\[
(I-\tau K)H^*=S^*+A^*-\tau KS^*-\tau KA^*,\qquad
H^{*\top}(I-\tau K)^\top=S^*-A^*+\tau S^*K-\tau A^*K,
\]
and summing the two identities yields \eqref{eq:sym-generalK}.
\end{proof}

\begin{proposition}
\label{prop:LRSGA-Jac}
Let $w^*$ be a Nash equilibrium and let $H^*=DF(w^*)=S^*+A^*$ with $S^*=S^{*\top}\succ0$ and $A^*=-A^{*\top}$.
Let $\alpha_k$ be the LRSGA skew correction at step $k$ (see (\ref{alphak})) and set $\Delta_k:=\alpha_k-A^*$ and $G_k:=(I-\tau\alpha_k)H^*$.
Consider 
$
T^{\mathrm{LRSGA}}_{k}(w)=w-\eta\,(I-\tau\alpha_k)\,F(w)$, whose Jacobian is $DT^{\mathrm{LRSGA}}_{k}(w^*)= I- \eta (I-\tau \alpha_k)H^*= I- \eta G_k $, and 
assume that 
\begin{equation}
\label{eq:tau-comm}
0<\tau<\frac{2\,\lambda_{\min}(S^*)}
{\ \|S^*A^*-A^*S^*\|_2 \;+\; \|S^*\Delta_k-\Delta_k S^*\|_2 \;+\; \|\Delta_k A^*+A^*\Delta_k\|_2\ },
\end{equation}

then, there exists $\bar\eta_k>0$ such that, for every $0<\eta<\bar\eta_k$, it holds
\(
\rho\, \!\big(D(T^{\mathrm{LRSGA}}_{k}(w^*))\big)\;<\;1.
\)
\end{proposition}

\begin{proof}
By Lemma~\ref{lem:sym-generalK} with $K=\alpha_k=A^*+\Delta_k$,
\[
\mathrm{sym}(G_k)=S^*+\tfrac{\tau}{2}(S^*A^*-A^*S^*)-\tau A^{*2}
+\tfrac{\tau}{2}\big(S^*\Delta_k-\Delta_kS^*-\Delta_kA^*-A^*\Delta_k\big).
\]
Since $-\,\tau (A^*)^{2}\succeq0$ we obtain
\[
\mathrm{sym}(G_k)\ \succeq\ c_k I,
\quad\text{with}\quad
c_k:=\lambda_{\min}(S^*)-\frac{\tau}{2}\Big(\|S^*A^*-A^*S^*\|_2
+\|S^*\Delta_k-\Delta_k S^*\|_2+\|\Delta_k A^*+A^*\Delta_k\|_2\Big).
\]
The condition \eqref{eq:tau-comm} implies $c_k>0$. 
If $\lambda$ is an eigenvalue of $G_k$, then
\[
|1-\eta\lambda|^2 \;\le\; 1-2\eta\,\Re\lambda+\eta^2|\lambda|^2
\;\le\; 1-2\eta c_k+\eta^2\|G_k\|_2^2.
\]
Setting
\[
\bar\eta_k:=\frac{2c_k}{\|G_k\|_2^2}>0,
\]
we obtain $|1-\eta\lambda|<1$ for every $0<\eta<\bar\eta_k$, hence $\rho(I-\eta G_k)<1$ and therefore
$\rho(DT^{\mathrm{LRSGA}}_{k}(w^*))<1$.
\end{proof}

As proved in \cite{lrsga}, letting $w^* \in \Omega$ be a NE point and assuming there is a $\delta > 0$ such that $\|{\mu_k -  D(\partial_{x} f(w^*))}\|_2 \leq \delta$ and $\|{\nu_k - D(\partial_{y} g(w^*))}\|_2 \leq \delta$ for $\mu_k$ and $\nu_k$ (see equations (\ref{eq:mu-nu-interif}) and (\ref{eq:mu-nu-interig})), it holds that the matrix $\alpha_k$, defined in (\ref{alphak}), satisfies the inequality $\|{\Delta_k}\|_2 \leq \delta$ with $\Delta_k:= \alpha_k - A^*$. Exploiting this result, the $\tau$-bound per iteration in Proposition \ref{prop:LRSGA-Jac} can be reformulated in global terms as follows.

\begin{corollary}
    Let $w^* \in \Omega$ and assume there is a $\delta > 0$ such that it holds $\|{\mu_k -  D(\partial_{x} f(w^*))}\|_2 \leq \delta$ and $\|{\nu_k - D(\partial_{y} g(w^*))}\|_2 \leq \delta$. Considering again $T^{\mathrm{LRSGA}}_{k}(w)=w-\eta\,(I-\tau\alpha_k)\,F(w)$, if 
    \[
    0<\tau<\frac{2\lambda_{\min}(S^*)}{\|S^*A^*-A^*S^*\|_2+2(\|S^*\|_2+\|A^*\|_2)\,\delta },
    \]
    then there exists $\bar\eta_k>0$ such that, for every $0<\eta<\bar\eta_k$, it holds 
    \(
    \rho\,\big(DT^{\mathrm{LRSGA}}_{k}(w^*)\big)\;<\;1.
    \)
\end{corollary}

\begin{proof}
Simply by calculus we have 
\[
\|S^*\Delta_k-\Delta_kS^*\|_2\le 2\|S^*\|_2\|\Delta_k\|_2,\qquad
\|\Delta_kA^*+A^*\Delta_k\|_2\le 2\|A^*\|_2\|\Delta_k\|_2, \qquad 
\]
By hypothesis we can apply \cite{lrsga} result about $\Delta_k$ obtaining:
\[
\|S^*\Delta_k-\Delta_kS^*\|_2\le 2\|S^*\|_2\delta,\qquad
\|\Delta_kA^*+A^*\Delta_k\|_2\le 2\|A^*\|_2\delta,
\]
and then the $\tau$–condition in
Proposition~\ref{prop:LRSGA-Jac} is implied by the global condition 
    \[
    0<\tau<\frac{2\lambda_{\min}(S^*)}{\|S^*A^*-A^*S^*\|_2+2(\|S^*\|_2+\|A^*\|_2)\,\delta }
    \]
\end{proof}

We remark that these results specialize to SGA in the limit $\Delta_k\!\to\!0$.
In particular, setting $\Delta_k=0$ yields
\[
\rho\,\!\bigl(DT^{\mathrm{SGA}}_{\star}(w^*)\bigr)<1
\quad\text{for}\quad
T^{\mathrm{SGA}}_{\star}(w):=w-\eta\,(I-\tau A^*)\,F(w),
\]
under the same restriction on $\tau$ and a corresponding choice of $\eta$.

These stability estimates, stated for the full LRSGA skew correction matrix $\alpha_k$, will be used in Section~5 to prove local convergence of LM--LRSGA, since the limited-memory skew correction matrix falls within the perturbation framework of Proposition~\ref{prop:LRSGA-Jac} (via Corollary~3.3).

\section{Limited Memory LRSGA }
\label{sec:lmlrsga}
We begin this section by introducing a further 'simplification' for LRSGA before address the problem the storage difficulties.
So, instead of working directly on LRSGA, we will work on this revised version.

It is algebraically and computationally equivalent (at the level of per–iteration arithmetic) to (i) form \(\mu_k\) and \(\nu_k\) and then extract only the mixed blocks \([\mu_k]_2\) and \([\nu_k]_1\), or (ii) update directly the mixed blocks. Indeed, starting from (\ref{eq:mu-nu-interif} and \ref{eq:mu-nu-interig}), by partitioning \(\mu_k = \bigl([\mu_k]_1 \,\big|\, [\mu_k]_2\bigr)\) and \(\nu_k = \bigl([\nu_k]_1 \,\big|\, [\nu_k]_2\bigr)\), one readily verifies
\begin{align} \label{mu-nu}
[\mu_{k+1}]_2 &= [\mu_{k}]_2  + \frac{\bigl(\partial_x f_{k+1} - \partial_x f_k -[\mu_k]_1(x_{k+1}-x_k)- [\mu_k]_2 (y_{k+1} - y_k)\bigr)(y_{k+1} - y_k)^\top}{(w_{k+1} - w_k)^\top(w_{k+1} - w_k) }  \\
[\nu_{k+1}]_1  &= [\nu_k]_1 + \frac{\bigl(\partial_y g_{k+1} - \partial_y g_k - [\nu_{k}]_2(y_{k+1}-y_k)- [\nu_k]_1  (x_{k+1} - x_k)\bigr)(x_{k+1} - x_k)^\top}{(w_{k+1} - w_k)^\top(w_{k+1} - w_k)}
\end{align}
where \( f_k :=  f(w_k)\) and \( g_k :=  g(w_k)\). 
%Note that
%\(
%(w_{k+1}-w_k)^\top(w_{k+1}-w_k)=\|x_{k+1}-x_k\|_2^2+\|y_{k+1}-y_k\|_2^2.
%\)
Thus, in either formulation the same quantities—namely the products \([\mu_k]_1(x_{k+1}-x_k)\) and \([\nu_k]_2(y_{k+1}-y_k)\)—are required, and the per–step computational cost is comparable. 
\\
In a neighborhood of a stable Nash equilibrium $w^{*}$, the second–order necessary
conditions imply that the own–player Hessian blocks are positive semidefinite,
\[
\partial^2_{xx} f(w^{*}) \succeq 0, 
\qquad 
\partial^2_{yy} g(w^{*}) \succeq 0.
\]
Locally around $w^{*}$, we assume that the diagonal blocks are dominated by their
coercive (positive) component, eventually corresponding to a stabilizing quadratic behavior of the form
\[
 f(x,y) + \tfrac{\lambda_x}{2}\|x-x^{*}\|^2,
\qquad
 g(x,y) + \tfrac{\lambda_y}{2}\|y-y^{*}\|^2,
\]
with $\lambda_x,\lambda_y > 0$.
Such  contributions do not alter rotational
effects in the local dynamics.

Since our objective is to eliminate the rotational component of the game dynamics,
we deliberately neglect the diagonal Hessian blocks and retain only the off–diagonal
coupling, which is responsible for non-potential behavior.
Accordingly, we approximate the local Jacobian $H(w_k)$ by its skew–symmetric part,
while assuming that the diagonal contribution is implicitly regularized by the
underlying positive curvature.

To this end, we introduce constant parameters $\varepsilon_x,\varepsilon_y \ge 0$,
interpreted as dominant positive eigenvalues of the diagonal blocks, and rewrite the
secant conditions for the off–diagonal approximations as
\begin{equation}\label{iterazconeps}
\begin{aligned}
M_{k+1} &=
M_k + 
\frac{
\big(
\partial_x f_{k+1} - \partial_x f_k
- \varepsilon_x (x_{k+1}-x_k)
- M_k (y_{k+1}-y_k)
\big)(y_{k+1}-y_k)^\top
}{
(w_{k+1}-w_k)^\top(w_{k+1}-w_k)
},
\\[6pt]
N_{k+1} &=
N_k + 
\frac{
\big(
\partial_y g_{k+1} - \partial_y g_k
- \varepsilon_y (y_{k+1}-y_k)
- N_k (x_{k+1}-x_k)
\big)(x_{k+1}-x_k)^\top
}{
(w_{k+1}-w_k)^\top(w_{k+1}-w_k)
}.
\end{aligned}
\end{equation}

The terms $\varepsilon_x (x_{k+1}-x_k)$ and $\varepsilon_y (y_{k+1}-y_k)$ model the
dominant coercive effect of the diagonal blocks and act as a built–in
Levenberg–Marquardt–type regularization, not affecting the rotational structure of the dynamics.

%%%%%%%%%%%%%%%%%%%%%%%%%%%%%%%%%%%%%%%%%%%%%%%%%%%%%%%%%

\begin{remark}
The introduction of $\varepsilon_x,\varepsilon_y$ in \eqref{iterazconeps} should be understood as a surrogate for the
diagonal (own--player) blocks that are implicitly present in the original secant relations.
Indeed, when one derives the mixed-block recursions induced by the full least--change updates, the formulas
contain the diagonal-block \emph{actions} $[\mu_k]_1(x_{k+1}-x_k)$ and $[\nu_k]_2(y_{k+1}-y_k)$.
The replacements
\[
[\mu_k]_1(x_{k+1}-x_k)\approx \varepsilon_x(x_{k+1}-x_k),
\qquad
[\nu_k]_2(y_{k+1}-y_k)\approx \varepsilon_y(y_{k+1}-y_k),
\]
provide the minimal block-level closure that preserves a secant-type update for the \emph{mixed} blocks while avoiding
the explicit formation and storage of the full matrices $\mu_k\in\mathbb{R}^{m\times(m+n)}$ and
$\nu_k\in\mathbb{R}^{n\times(m+n)}$. In this way one can update and store directly
$M_k\in\mathbb{R}^{m\times n}$ and $N_k\in\mathbb{R}^{n\times m}$ (instead of computing $\mu_k,\nu_k$ and then
extracting only $[\mu_k]_2$ and $[\nu_k]_1$), while still accounting for the dominant coercive effect of the diagonal
blocks through $\varepsilon_x,\varepsilon_y$.
Finally, this approximation does \emph{not} change the equilibrium conditions of the iteration: since $\alpha_k$ is skew-symmetric, $I-\tau\alpha_k$ is nonsingular for every $\tau>0$, hence
$ (I-\tau\alpha_k)F(w_k)=0 \iff F(w_k)=0$. 
%Thus, the $\varepsilon$-terms only affect the quality of the mixed-block approximation (and consequently the constants in the bounds).
\end{remark}

%%%%%%%%%%%%%%%%%%%%%%%%%%%%%%%%%%%%%%%%%%%%%%%%%%%%%%%%%%%%%%
As a result, the correction matrix adopted in this slightly revised LRSGA is purely skew–symmetric,
\begin{equation}\label{alphaknew}
\overline{\alpha}_k =
\begin{pmatrix}
0 & \tfrac12(M_k - N_k^{\top}) \\[6pt]
\tfrac12(N_k - M_k^{\top}) & 0
\end{pmatrix},
\end{equation}
which isolates and counteracts the rotational component of the game Jacobian,
while the stabilizing diagonal curvature is implicitly absorbed into the step
regularization.

In any case, as highlighted, LRSGA requires storing the matrices \(M_k\) and \(N_k\) to approximate the mixed second derivatives, which can be memory–intensive for high–dimensional models. To address this, we propose a \emph{Limited-Memory LRSGA} variant that preserves LRSGA’s theoretical robustness while substantially reducing memory usage.
Our strategy follows the \emph{two-loop recursion} method, successfully employed by L-BFGS \cite{nocedal1980updating,liu1989limited}. The key idea is to avoid storing the full updated matrices \(M_k\) and \(N_k\); instead, we retain only the most recent \(\ell\) pairs among the four vectors \(\{s_k^x, s_k^y, y_k^f, y_k^g\}\), defined as
\[
s_k^x = x_{k+1} - x_{k}, \qquad
s_k^y = y_{k+1}  - y_{k},
\]
\[
y_k^f = \partial_x f_{k+1}  - \partial_x f_{k} - \varepsilon_x\, s_k^x, \qquad
y_k^g = \partial_y g_{k+1}  - \partial_y g_{k} - \varepsilon_y\, s_k^y.
\]
Adopting this notation and denoting \(s_k^{w} := (s_k^{x}, s_k^{y})\), LRSGA can be written as follows:

\begin{algorithm}[H]
\caption{LRSGA (Low-Rank Symplectic Gradient Adjustment)}
\label{alg:LRSGA}
\begin{algorithmic}[1]
\State \textbf{Phase 1: Gradient Computation and Update}
\For{\(k = 1,2,\dots\)}
    \State \(M_{k} \gets M_{k-1} + \frac{\bigl(y_{k-1}^f - M_{k-1}\,s_{k-1}^y\bigr)}{(s_{k-1}^w)^{\top} s_{k-1}^w}\,(s_{k-1}^y)^{\top}\)
    \State \(N_{k} \gets N_{k-1} + \frac{\bigl(y_{k-1}^g - N_{k-1}\,s_{k-1}^x\bigr)}{(s_{k-1}^w)^{\top} s_{k-1}^w}\,(s_{k-1}^x)^{\top}\)
    
    \State \(x_{k+1} \gets x_k - \eta \Bigl( \partial_x f_k - \frac{\tau}{2}\bigl(M_k\, - N_k^{\top}\,\bigr) \partial_y g_k\Bigr)\)
    \State \(y_{k+1} \gets y_k - \eta \Bigl( \partial_y g_k - \frac{\tau}{2}\bigl(N_k\, - M_k^{\top}\,\bigr)\partial_x f_k \Bigr)\)
\EndFor
\State \textbf{Output:} Updated \(x_k\) and \(y_k\) until convergence.
\end{algorithmic}
\end{algorithm}

\subsection[Limited-memory recursive updates for Mk and Nk]%
{Limited-memory recursive updates for $M_k$ and $N_k$}

The updates for the matrices $M_k$ and $N_k$ can be written in a limited-memory, iterative form by storing only the last $\ell$ pairs $(s_t,y_t)$ and using the auxiliary matrices
\[
\widetilde V_k := I - p_k\, s_k^y (s_k^y)^\top, 
\qquad 
V_k := I - p_k\, s_k^x (s_k^x)^\top,
\]
with
\[
p_k := \frac{1}{(s_k^w)^\top s_k^w}.
\]
Consequently, iterations \ref{iterazconeps} can be rewritten as:
\[
M_{k} = M_{k-1} \widetilde V_{k-1}+ p_{k-1} y_{k-1}^f(s_{k-1}^y)^{\top} \qquad \text{and } \quad N_{k} = N_{k-1} {V_{k-1}}+ p_{k-1} y_{k-1}^g(s_{k-1}^x)^{\top}
\]
By considering only the last $\ell$ pairs $(s_t,y_t)$, these became: 
\begin{equation} 
\begin{aligned}
M_{k}^{(\ell)} &= M_{k-\ell} (\widetilde V_{k-\ell}\dots \widetilde V_{k-1}) + p_{k-\ell} y_{k-\ell}^f (s_{k-\ell}^y)^{\top}(\widetilde V_{k-\ell+1}\dots \widetilde V_{k-1}) + \dots \\
&\quad+ p_{k-2} y_{k-2}^f (s_{k-2}^y)^{\top} \widetilde V_{k-1} + p_{k-1} y_{k-1}^f(s_{k-1}^y)^{\top},
\end{aligned}
\end{equation}

\begin{equation}
\begin{aligned}
N_{k}^{(\ell)} &= N_{k-\ell} ({V}_{k-\ell}\dots {V}_{k-1}) + p_{k-\ell} y_{k-\ell}^g (s_{k-\ell}^x)^{\top}({V}_{k-\ell+1}\dots {V}_{k-1}) + \dots \\
&\quad+ p_{k-2} y_{k-2}^g (s_{k-2}^x)^{\top} {V}_{k-1} + p_{k-1} y_{k-1}^g(s_{k-1}^x)^{\top}.
\end{aligned}
\end{equation}

In practice, we avoid storing the base matrices $M_{k-\ell}$ and $N_{k-\ell}$, by replacing them with respectively the rank-one initializations: 
\begin{equation}
    \label{eq:inizializzaioni}
H_0 ^k := p_{k-\ell} \, y_{k- \ell}^{\,h}\,(s_{k -\ell }^{\,j})^\top,
\end{equation}
where \(j\in\{x,y\}\) an $h \in \{f,g\}$ matches the recursion (i.e., \(j=y\) and $h=f$ for the \(M_k\)-recursion, and \(j=x\), $h=g$ for the \(N_k\)-recursion.

Then we can redefine the limited-memory approximations $M_k^{(\ell)}$ and $N_k^{(\ell)}$ as

\begin{equation} \label{eq:MKesattaNonfull}
\begin{aligned}
M_{k}^{(\ell)} &=H_0 ^k (\widetilde V_{k-\ell}\dots \widetilde V_{k-1}) + p_{k-\ell} y_{k-\ell}^f (s_{k-\ell}^y)^{\top}(\widetilde V_{k-\ell+1}\dots \widetilde V_{k-1}) + \dots \\
&\quad+ p_{k-2} y_{k-2}^f (s_{k-2}^y)^{\top} \widetilde V_{k-1} + p_{k-1} y_{k-1}^f(s_{k-1}^y)^{\top},
\end{aligned}
\end{equation}

\begin{equation}\label{eq:NKesattaNonfull}
\begin{aligned}
N_{k}^{(\ell)} &= H_0 ^k ({V}_{k-\ell}\dots {V}_{k-1}) + p_{k-\ell} y_{k-\ell}^g (s_{k-\ell}^x)^{\top}({V}_{k-\ell+1}\dots {V}_{k-1}) + \dots \\
&\quad+ p_{k-2} y_{k-2}^g (s_{k-2}^x)^{\top} {V}_{k-1} + p_{k-1} y_{k-1}^g(s_{k-1}^x)^{\top}.
\end{aligned}
\end{equation}

By symmetry of the projectors ($\widetilde V_t=\widetilde V_t^{\top}$ and $V_t=V_t^{\top}$), the compact transposed forms read:
\[
(M_{k}^{(\ell)})^{\top} = \Bigl( \widetilde{V}_{k-1} \cdots \widetilde{V}_{k-\ell} \Bigr) (H_0 ^k)^\top + p_{k-\ell} \Bigl( \widetilde{V}_{k-1} \cdots \widetilde{V}_{k-\ell+1} \Bigr) s_{k-\ell}^y\, (y_{k-\ell}^f)^{\top} + \cdots + p_{k-1}\, s_{k-1}^y\, (y_{k-1}^f)^{\top},
\]
and, analogously,
\[
(N_{k}^{(\ell)})^{\top} = \Bigl( V_{k-1} \cdots V_{k-\ell} \Bigr) (H_0 ^k)^\top + p_{k-\ell} \Bigl( V_{k-1} \cdots V_{k-\ell+1} \Bigr) s_{k-\ell}^x\, (y_{k-\ell}^g)^{\top} + \cdots + p_{k-1}\, s_{k-1}^x\, (y_{k-1}^g)^{\top}.
\]

These representations enable the computation of the four products required by LRSGA—namely
\(M_k\,\partial_y g(w_k),\; N_k^{\top}\,\partial_y g(w_k),\; N_k\,\partial_x f(w_k),\; M_k^{\top}\,\partial_x f(w_k)\)—while substantially reducing the memory footprint.

\begin{remark}[Choice of the rank-one base $H_0^k$]
In the limited-memory implementation we maintain a buffer of length $\ell$ containing the most recent curvature pairs, i.e., indices $t \in \{k-\ell,\dots,k-1\}$ once the buffer is full. Accordingly, in \eqref{eq:inizializzaioni} we choose
the rank-one base matrix $H_0^k$ using the \emph{oldest pair within the current window}, namely $t=k-\ell$. With this choice, the oldest curvature pair $(s_{k-\ell}^j,y_{k-\ell}^h)$ appears twice in the explicit expansions \eqref{eq:MKesattaNonfull}--\eqref{eq:NKesattaNonfull}: once through the base term $H_0^k(\widetilde V_{k-\ell}\cdots \widetilde V_{k-1})$ (respectively $H_0^k(V_{k-\ell}\cdots V_{k-1})$), and once through the first rank-one term in the truncated sum. As an alternative, one may initialize the base using the curvature pair immediately \emph{preceding} the memory window,
\[
\hat{H}_0^k := p_{k-\ell-1}\,y_{k-\ell-1}^{\,h}\,(s_{k-\ell-1}^{\,j})^\top .
\]
In that case one recovers the telescoping identity
\(
M_k - M_k^{(\ell)} = M_{k-\ell-1}\,(\widetilde V_{k-\ell-1}\cdots \widetilde V_{k-1})\), \(
N_k - N_k^{(\ell)} = N_{k-\ell-1}\,(V_{k-\ell-1}\cdots V_{k-1}),
\)
and therefore the sharper bounds $\|M_k-M_k^{(\ell)}\|_2 \le \|M_{k-\ell-1}\|_2$ and $\|N_k-N_k^{(\ell)}\|_2 \le \|N_{k-\ell-1}\|_2$ (using $\|\widetilde V_t\|_2\le 1$ and $\|V_t\|_2\le 1$).
From an implementation standpoint, this alternative requires computing  $\hat{H}_0^k$ \emph{before} the $(k-\ell-1)$-th curvature pair is discarded, i.e., prior to updating the history lists (i.e., before the pop operation).
In this paper (and in our experiments) we adopt the implementation-consistent choice $H_0^k$ \eqref{eq:inizializzaioni} with
index $k-\ell$, so that no curvature pair beyond the last $\ell$ needs to be accessed or stored.
All subsequent norm bounds and local convergence arguments only require a rank-one base of uniformly bounded norm,
hence they apply unchanged to either initialization.
\end{remark}

\subsection{Adapted two-loop recursion algorithms}

Classical L-BFGS uses the two–loop recursion to apply a \emph{square} (inverse) Hessian approximation $H_k\in\mathbb{R}^{d\times d}$ to a vector. 
In contrast, our setting requires multiplying \emph{rectangular} mixed-Jacobian approximation and their transpositions,
$M_k\in\mathbb{R}^{m\times n}$ and $N_k\in\mathbb{R}^{n\times m}$, by vectors of compatible size.
To this end, we introduced the two following two-loop recursion algorithms.

\subsubsection{Direct Two-loop recursion }

Given a generic vector $q$, we can compute $M_k q$ and  $N_k q$ via the Algorithm~\ref{alg:twoloop_direct} below. 

\begin{algorithm}[H]
\caption{TwoLoopRecursionDirect-K}\label{alg:twoloop_direct}
\begin{algorithmic}[1]
\State \textbf{Input:} $\{p_{k-1}, p_{k-2}, \dots , p_{k-\ell}\},\{s_{k-1}^j, s_{k-2}^j, \dots , s_{k-\ell}^j\}, \{y_{k-1}^h, y_{k-2}^h, \dots , y_{k-\ell}^h\}$, $j\in\{x,y\}$, $h\in \{f,g\}$, vector $q$, initial matrix $H_0^k$
\For {$i = k-1, k-2, \dots, k-\ell$}
\State $\alpha_i \leftarrow p_i(s_i^j)^{\top} q$
\State $q \leftarrow q - \alpha_i s_i^j$
\EndFor
\State $r \leftarrow H_0^k q$
\For {$i = k-\ell, k-\ell+1, \dots, k-1$}
\State $r \leftarrow r + y_i^h \alpha_i$
\EndFor
\State \textbf{Output:} $r$
\end{algorithmic}
\end{algorithm}

Applying Algorithm~\ref{alg:twoloop_direct} with the appropriate inputs yields the LRSGA products:
\[
\resizebox{\textwidth}{!}{$
\boxed{\;
M_k\,\partial_y g(w_k)
\;=\;
\textbf{TwoLoopRecursionDirect-K}\bigl(\{p_t,s_t^y,y_t^f\}_{t=k-\ell}^{k-1},\ q=\partial_y g(w_k),\ H_0^k\bigr),
\;} $}
\]
\[
\resizebox{\textwidth}{!}{$
\boxed{\;
N_k\,\partial_x f(w_k)
\;=\;
\textbf{TwoLoopRecursionDirect-K}\bigl(\{p_t,s_t^x,y_t^g\}_{t=k-\ell}^{k-1},\ q=\partial_x f(w_k),\ H_0^k \bigr).
\;} $}
\]

\subsubsection{Transpose Two-loop recursion }

Given a vector $q$, to compute products $M_k^{\top} q$ and $N_k^{\top} q$), we use the following limited–memory transpose recursions that mirror the direct case. 

\begin{algorithm}[H]
\caption{TwoLoopRecursionTranspose-K
}\label{alg:twoloop_transpose}
\begin{algorithmic}[1]
\State \textbf{Input:} $\{p_{k-1}, p_{k-2}, \dots , p_{k-\ell}\},\{s_{k-1}^j, s_{k-2}^j, \dots , s_{k-\ell}^j\}, \{y_{k-1}^h, y_{k-2}^h, \dots , y_{k-\ell}^h\}$, $j\in\{x,y\}$, $h\in \{f,g\}$,  vector $q$, initial  matrix  $H_0^k$
\For {$i = k-1, k-2, \dots, k-\ell$}
\State $\alpha_i \leftarrow p_i(y_i^h)^{\top} q$
\EndFor
\State $r \leftarrow (H_0^k)^{\top} q$
\For {$i = k-\ell, k-\ell+1, \dots, k-1$}
\State $\beta_i \leftarrow p_i(s_i^j)^{\top} r$
\State $r \leftarrow r + s_i^j(\alpha_i-\beta_i)$
\EndFor
\State \textbf{Output:} $r$
\end{algorithmic}
\end{algorithm}

In particular, we apply Algorithm~\ref{alg:twoloop_transpose}, as follows:

\[
\resizebox{\textwidth}{!}{$
\boxed{\;
M_k^{\top}\,\partial_x f(w_k)
\;=\;
\textbf{TwoLoopRecursionTranspose-K}\bigl(\{p_t,s_t^y,y_t^f\}_{t=k-\ell}^{k-1},\ q=\partial_x f(w_k),\ H_0^k \bigr),
\;} $}
\]
\[
\resizebox{\textwidth}{!}{$
\boxed{\;
N_k^{\top}\,\partial_y g(w_k)
\;=\;
\textbf{TwoLoopRecursionTranspose-K}\bigl(\{p_t,s_t^x,y_t^g\}_{t=k-\ell}^{k-1},\ q=\partial_y g(w_k),\ H_0^k)^{\top}\bigr).
\;} $}
\]

%Per quanto riguarda la scelta di $M_k^0$ and  $N_k^0$ abbiamo optato per:

\begin{comment}
\begin{remark}\label{rem:telesc}
Fix $\ell\ge 1$ and consider the $\ell$–step truncated update $M_k^{(\ell)}$ by replacing $M_{k-\ell}$ with a base matrix
$H_0^k$ we obtain
\begin{align*}
M_k^{(\ell)}
&= H_0^k\,(\widetilde V_{k-\ell}\cdots \widetilde V_{k-1})
\;+\; p_{k-\ell}\,y_{k-\ell}^{f}\,(s_{k-\ell}^{y})^{\!\top}\,(\widetilde V_{k-\ell+1}\cdots \widetilde V_{k-1})
\;+\; \cdots \\
&\;+\; p_{k-2}\,y_{k-2}^{f}\,(s_{k-2}^{y})^{\!\top}\,\widetilde V_{k-1}
\;+\; p_{k-1}\,y_{k-1}^{f}\,(s_{k-1}^{y})^{\!\top}.
\end{align*}
Then the the difference between the full $M_k$ and the $\ell$-step truncated one is
\[
M_k-M_k^{(\ell)}=(M_{k-\ell}-H_0^k)\,(\widetilde V_{k-\ell}\cdots \widetilde V_{k-1}).
\]
In particular, with the chosen rank–one base
\(
H_0^k:=p_{k-\ell}\,y_{k-\ell}^f (s_{k-\ell}^y)^\top,
\)
 and since for the projectors $\|\widetilde V_t\|_2\le 1$, we obtain 
\[
\|M_k-M_k^{(\ell)}\|_2 \;\le\; \|M_{k-\ell}-H_0^k\|_2 \le \|M_{k-\ell}\|_2+\|H_0^k\|_2.
\]

The same identities and bounds hold for $N_k$ with $V_t$.
\smallskip

\end{remark}
\end{comment}

We remark that in both differences $\partial_x f(w_{k+1})-\partial_x f(w_k)$ and $\partial_y g(w_{k+1})-\partial_y g(w_k)$
relate two consecutive iterates and this implies that each term is evaluated on the \emph{current} mini–batch $B$
(e.g., \(\partial_x f_{B_{k+1}}(w_{k+1})-\partial_x f_{B_k}(w_k)\)).
We emphasize that this fact degrades the accuracy of the
mixed–second–derivative approximations used by our secant updates and so a stochastic implementation is required.

In line with stochastic quasi–Newton practice \cite{guo2023overview}, we therefore adopt
\emph{gradient displacement} (also called common/consistent–batch ): both gradients in each difference are computed on
the \emph{same} mini–batch, i.e.,
\[
\partial_x f_{B_k}(w_{k+1})-\partial_x f_{B_k}(w_k),
\qquad
\partial_y g_{B_k}(w_{k+1})-\partial_y g_{B_k}(w_k),
\]
which reduces variance and preserves secant consistency for the mixed blocks.
As an alternative, one may enforce overlapping mini–batches with
\(I_k:=B_k\cap B_{k+1}\neq\emptyset\) and form the differences on the
intersection,
\[
\partial_x f_{I_k}(w_{k+1})-\partial_x f_{I_k}(w_k),
\qquad
\partial_y g_{I_k}(w_{k+1})-\partial_y g_{I_k}(w_k).
\]

Since common/consistent-batch gradient displacement can increase implementation complexity (and may require additional gradient evaluations at $w_{k+1}$ on the previous mini-batch $B_k$), we propose an additional variant, \textbf{LM-LRSGAEma}, which mitigates mini-batch–induced noise without manipulating the batches. Specifically, we maintain exponential moving averages of $y_k^{\,f}$ and $y_k^{\,g}$ and use them in the two–loop recursions:
\[
\tilde y_k^{\,f} = \beta\,\tilde y_{k-1}^{\,f} + (1-\beta)\,y_k^{\,f}, 
\qquad
\tilde y_k^{\,g} = \beta\,\tilde y_{k-1}^{\,g} + (1-\beta)\,y_k^{\,g},
\quad \beta\in[0,1),
\]
replacing $y_k^{\,f},y_k^{\,g}$ with $\tilde y_k^{\,f},\tilde y_k^{\,g}$ in the update formulas. This EMA version seems to reduces variance, and improves stability without altering the underlying dynamics.

In conclusion, the proposed approach computes LRSGA updates with a markedly reduced memory footprint, thereby substantially broadening its applicability to high–dimensional models.

\section{On the local convergence of LM-LRSGA}
\label{sec:convergence-lmlrsga}

In this section we prove a local linear convergence guarantee for the limited--memory variant introduced in Section~\ref{sec:lmlrsga}. Specifically, for any fixed history length~$\ell$, and for sufficiently small stepsizes
$\eta$ and $\tau$, the LM-LRSGA iterates converge linearly to a stable Nash equilibrium, with constants depending on~$\ell$.

%In this section we show that the limited–memory variant introduced in Section~\ref{sec:lmlrsga} enjoys the same type of local linear convergence guarantee as the full LRSGA scheme, provided that the history length~$\ell$ is fixed and the stepsizes $\eta$ and $\tau$ are sufficiently small.

Throughout the section we work on a convex open set
$\Omega\subset\mathbb{R}^{m+n}$ and assume
$f,g\in C^3(\Omega,\mathbb{R})$, so that the game gradient
$F(w) = (\partial_x f(w),\partial_y g(w))^\top$ is of class $C^2$.

Recall that the LM-LRSGA iteration can be written as
\begin{equation}\label{eq:LM-iteration}
    w_{k+1} = w_k - \eta\bigl(I - \tau\,\alpha_k^{(\ell)}\bigr)F(w_k),
\end{equation}
where $\eta>0$ and $\tau>0$ are fixed stepsizes and $\alpha_k^{(\ell)}$ is the limited–memory approximation of the antisymmetric
component of the game Hessian
\begin{equation}\label{eq:alpha-lm}
    \alpha_k^{(\ell)}
    :=
    \begin{pmatrix}
        0 & \tfrac12\bigl(M_k^{(\ell)} - (N_k^{(\ell)})^\top\bigr) \\[2mm]
        \tfrac12\bigl(N_k^{(\ell)} - (M_k^{(\ell)})^\top\bigr) & 0
    \end{pmatrix},
\end{equation} obtained from the truncated mixed–blocks
$M_k^{(\ell)}$ and $N_k^{(\ell)}$ (see \ref{eq:MKesattaNonfull} and \ref{eq:NKesattaNonfull}) 
with curvature pairs
$(s_k^x,s_k^y,y_k^f,y_k^g)$ and fixed positive parameters
$\varepsilon_x,\varepsilon_y$.

We start by showing, in the following lemma, that the limited-memory mixed blocks  $M_k^{(\ell)}$ and $N_k^{(\ell)}$ are uniformly bounded.

\begin{lemma}\label{lem:lm-bounded-MN}
Assume that $\partial_x f$ and $\partial_y g$ are Lipschitz on $\Omega$ with
constants $L_x,L_y>0$.
%, and that the iterates $w_k$ generated by LM-LRSGA remain in~$\Omega$.
%Let $M_k^{(\ell)}$ and $N_k^{(\ell)}$ be defined by the limited–memory recursions \eqref{16}--\eqref{17}, .
Then, for every fixed history length $\ell\ge 1$ there exist constants
\[
    C_M^x(\ell) := (\ell+1)\bigl(L_x + |\varepsilon_x|\bigr),
    \qquad
    C_M^y(\ell) := (\ell+1)\bigl(L_y + |\varepsilon_y|\bigr),
\]
such that, for all $k$ with $s_k^w\neq 0$,
\[
    \|M_k^{(\ell)}\|_2 \le C_M^x(\ell),
    \qquad
    \|N_k^{(\ell)}\|_2 \le C_M^y(\ell),
\]
where $s_k^w := (s_k^x,s_k^y)$.
\end{lemma}

\begin{proof}
Fix $k\ge 1$. By Lipschitz continuity of $\partial_x f$ on~$\Omega$,
\[
    \|\partial_x f(w_{k+1} ) - \partial_x f(w_{k})\|_2
    \le L_x \|w_{k+1}  - w_{k}\|_2
    = L_x \|s_k^w\|_2.
\]
Using the definition of $y_k^f$ and the inequality
$\|s_k^x\|_2 \le \|s_k^w\|_2$ we obtain
\[
    \|y_k^f\|_2
    \le
    \|\partial_x f(w_{k+1} ) - \partial_x f(w_{k})\|_2
    + |\varepsilon_x|\,\|s_k^x\|_2
    \le
    \bigl(L_x + |\varepsilon_x|\bigr)\,\|s_k^w\|_2.
\]
Recalling that $p_k := 1/\|s_k^w\|_2^2$ and considering that $\|s_k^y\|_2\le\|s_k^w\|_2$, the rank–one term satisfies
\[
    \|p_k y_k^f (s_k^y)^\top\|_2
    \le
    p_k \|y_k^f\|_2 \,\|s_k^y\|_2
    \le
    \frac{\bigl(L_x+|\varepsilon_x|\bigr)\,
          \|s_k^w\|_2\,\|s_k^y\|_2}{\|s_k^w\|_2^2}
    \le
    L_x + |\varepsilon_x|.
\]
An identical argument applied to $\partial_y g$ and $y_k^g$ yields
\[
    \|p_k y_k^g (s_k^x)^\top\|_2
    \le
    L_y + |\varepsilon_y|.
\]

The limited–memory formulas (\ref{eq:MKesattaNonfull}) and (\ref{eq:NKesattaNonfull}) express
$M_k^{(\ell)}$ and $N_k^{(\ell)}$ as finite sums of at most $\ell+1$ such
rank–one terms, multiplied on the right by products of the matrices
\[
    \tilde V_t := I - p_t\,s_t^y (s_t^y)^\top,
    \qquad
    V_t := I - p_t\,s_t^x (s_t^x)^\top,
\]
for indices $t$ in the history window.
Since $\|s_t^x\|_2,\|s_t^y\|_2\le\|s_t^w\|_2$ and $p_t = \|s_t^w\|_2^{-2}$,
the eigenvalues of $\tilde V_t$ and $V_t$ lie in $[0,1]$, hence
$\|\tilde V_t\|_2\le 1$ and $\|V_t\|_2\le 1$ for all~$t$.
The initializations $H_0^k$ (see Equation\eqref{eq:inizializzaioni}) are also of the same form,
and therefore satisfy the same bounds.

Then, by the triangle inequality we
conclude that
\[
    \|M_k^{(\ell)}\|_2
    \le
    (\ell+1)\,\bigl(L_x + |\varepsilon_x|\bigr)
    = C_M^x(\ell),
\]
and analogously
\[
    \|N_k^{(\ell)}\|_2
    \le
    (\ell+1)\,\bigl(L_y + |\varepsilon_y|\bigr)
    = C_M^y(\ell),
\]
for every~$k$ with $s_k^w\neq 0$.
\end{proof}

By denoting
\(
    \Delta_k^{(\ell)} := \alpha_k^{(\ell)} - A^*, 
\)
the next corollary shows that $\Delta_k^{(\ell)}$ remains uniformly bounded along the iterations.

\begin{corollary}\label{cor:lm-bounded-alpha}
Under the assumptions of Lemma~\ref{lem:lm-bounded-MN}, let $w^*$ be a
stable Nash equilibrium and denote by $A^*$ the antisymmetric part of
the game Hessian at $w^*$.
Then there exists a constant $\delta^{(\ell)}>0$ 
% ,depending only on $\ell,L_x,L_y,\varepsilon_x,\varepsilon_y$ and $\|A^*\|_2$, 
such that
\[
    \|\alpha_k^{(\ell)} - A^*\|_2 \le \delta^{(\ell)}
    \qquad\text{for all }k.
\]
In particular, we may take
\[
    C_\alpha(\ell) :=
        \tfrac12\bigl(C_M^x(\ell)+C_M^y(\ell)\bigr),
    \qquad
    \delta^{(\ell)} := C_\alpha(\ell) + \|A^*\|_2 .
\]
\end{corollary}

\begin{proof}
As shown in~\cite{lrsga}, for any rectangular matrix $C$ the block–skew
matrix
\(
    B = \begin{psmallmatrix} 0 & C \\ -C^\top & 0 \end{psmallmatrix}
\)
satisfies $\|B\|_2 = \|C\|_2$.
Applying this to the definition \eqref{eq:alpha-lm} we obtain
\[
    \|\alpha_k^{(\ell)}\|_2
    = \tfrac12 \bigl\|M_k^{(\ell)} - (N_k^{(\ell)})^\top\bigr\|_2
    \le \tfrac12\bigl(\|M_k^{(\ell)}\|_2 + \|N_k^{(\ell)}\|_2\bigr)
    \le C_\alpha(\ell)
\]
for all~$k$, where $C_\alpha(\ell) :=
        \tfrac12\bigl(C_M^x(\ell)+C_M^y(\ell)\bigr)$.
Finally,
\[
    \|\Delta_k^{(\ell)}\|_2 := \|\alpha_k^{(\ell)} - A^*\|_2
    \le \|\alpha_k^{(\ell)}\|_2 + \|A^*\|_2
    \le C_\alpha(\ell) + \|A^*\|_2
    = \delta^{(\ell)},
\]
which proves the claim.
\end{proof}

We finally combine Corollary~\ref{cor:lm-bounded-alpha} with the spectral condition of Section~3 to obtain a local convergence result for LM-LRSGA as claimed in the following theorem.

Letting $w^*\in\Omega$ be stable Nash equilibrium, for each $k$ and fixed history length $\ell$, we can define the $k$-frozen LM-LRSGA map
\[
    T_{k,(\ell)}^{LM-LRSGA}(w) := w - \eta\,(I - \tau \alpha_k^{(\ell)})\,F(w),
\]
so that $w_{k+1} =  T_{k,(\ell)}^{LM-LRSGA}(w_k)$.
Its Jacobian at the equilibrium $w^*$ is
\[
    D T_{k,(\ell)}^{LM-LRSGA}(w^*)
    = I - \eta G_k^{(\ell)},
    \qquad
    G_k^{(\ell)} := (I - \tau \alpha_k^{(\ell)})\,H^*.
\]

In particular, the following theorem shows that, suitable choices of $\tau$ and $\eta$
yield local linear convergence of the full LM-LRSGA iteration.

%%%%%%%%%%%%%%%%%%%%%%%%%%%%%%%%%%%%%%%%%%%%%%%%%%%%%%%%%%%%%%%%%%%%%%%%%%%%%%%%%%%%%%%%%%%%%%%%%%%%%%%%%%%%%%%%%%%%%%%%%%%%%%%%%%%%%%%%%%%%%%%%%%%
\begin{theorem}
\label{thm:lm-lrsga-local-conv}
Let $\Omega\subset\mathbb{R}^{m+n}$ be a convex open set and $f,g\in C^3(\Omega,\mathbb{R})$.
Assume that $\partial_x f$ and $\partial_y g$ are Lipschitz on $\Omega$ with
constants $L_x,L_y>0$, and let $w^*\in\Omega$ be a stable Nash equilibrium.
Then, for any fixed history length $\ell\ge 1$ there exist constants
\[
    \bar\tau(\ell) > 0,\qquad
    \bar\eta(\ell) > 0,\qquad
    r(\ell) > 0,
\]
such that the following holds.
If $0 < \tau < \bar\tau(\ell)$, $0 < \eta < \bar\eta(\ell)$ and
$w_0\in B_{r(\ell)}(w^*)$, then:
\begin{enumerate}
    \item the LM-LRSGA iterates $w_{k+1} = T_{k,(\ell)}^{LM-LRSGA}(w_k)$ are well-defined
          and remain in $B_{r(\ell)}(w^*)$ for all $k$;
    \item the sequence $(w_k)$ converges linearly to $w^*$, i.e.\ there exists
          a constant $q\in(0,1)$ such that
          \[
              \|w_k - w^*\|_2 \le q^k \,\|w_0 - w^*\|_2
              \qquad\text{for all }k\ge 0.
          \]
\end{enumerate}
\end{theorem}

\begin{proof}
By Corollary~\ref{cor:lm-bounded-alpha}, for each fixed $\ell$ there exists
$\delta^{(\ell)}>0$ such that
$\|\alpha_k^{(\ell)} - A^*\|_2 \le \delta^{(\ell)}$ for all $k$.
Applying Corollary~3.3 with $\Delta_k^{(\ell)} := \alpha_k^{(\ell)} - A^*$ we
obtain a bound
\[
    0 < \tau < \bar\tau(\ell)
    :=
    \frac{2\lambda_{\min}(S^*)}
         {\|S^*A^* - A^*S^*\|_2 +
          2(\|S^*\|_2 + \|A^*\|_2)\,\delta^{(\ell)}}
\]
such that, for every such~$\tau$, the symmetric part of $G_k^{(\ell)}$ satisfies
\(
\operatorname{sym}\bigl(G_k^{(\ell)}\bigr)
    \succeq c(\ell)\,I
    \quad\text{for all }k,
\)
for some constant $c(\ell)>0$ depending on $\ell$ and the spectral data of
$H^*$. 
Moreover, from Lemma~\ref{lem:lm-bounded-MN} and the definition of
$\alpha_k^{(\ell)}$ we also have a uniform bound
\[
    \|G_k^{(\ell)}\|_2
    \le \bigl\|I - \tau\,\alpha_k^{(\ell)}\bigr\|_2 \,\|H^*\|_2
    \le (1+ \tau C_\alpha (\ell)) \|H^*\|_2 =: G_{\max}(\ell)
    \qquad\text{for all }k,
\]
for a suitable constant $G_{\max}(\ell)>0$.

For any $z\in\mathbb{R}^{m+n}$,
\begin{align*}
    %\| D(T_{k,(\ell)}^{LMLRSGA}(w^*)) z \| =
    \|(I - \eta G_k^{(\ell)})z\|_2^2
    &= \|z\|_2^2
       - 2\eta\,\bigl\langle \operatorname{sym}(G_k^{(\ell)})z,z\bigr\rangle
       + \eta^2\|G_k^{(\ell)}z\|_2^2 \\
    &\le \bigl(1 - 2\eta c(\ell) + \eta^2 G_{\max}(\ell)^2\bigr)\,\|z\|_2^2 .
\end{align*}
Choosing
\[
    0 < \eta < \bar\eta(\ell)
    := \frac{2c(\ell)}{G_{\max}(\ell)^2},
\]
we obtain
\(
    1 - 2\eta c(\ell) + \eta^2 G_{\max}(\ell)^2 < 1
\),
and therefore a uniform bound
\[
    \|I - \eta G_k^{(\ell)}\|_2
    \le q_0(\ell) < 1
    \qquad\text{for all }k,
\]
for some $q_0(\ell)\in(0,1)$ depending on
$\ell,\tau,H^*$ and the constants above.

Since $f,g\in C^3(\Omega)$, the game gradient $F$ is of class $C^2$.
Thus its Jacobian $DF(w)$ is continuous, and in a sufficiently small
ball $B_{r_0}(w^*)\subset\Omega$ we can use the Taylor expansion
\[
    F(w) = F(w^*) + H^*(w-w^*) + R(w),
\]
with remainder satisfying
\[
    \|R(w)\|_2 \le C_R \|w-w^*\|_2^2,
    \qquad w\in B_{r_0}(w^*),
\]
for some constant $C_R>0$.

For $w_k\in B_{r_0}(w^*)$ we have $F(w^*)=0$ and
\[
    w_{k+1}-w^*
    = (w_k - w^*) - \eta \bigl(I - \tau\alpha_k^{(\ell)}\bigr) F(w_k)
    = (I - \eta G_k^{(\ell)})(w_k - w^*)
      - \eta\bigl(I - \tau\alpha_k^{(\ell)}\bigr)R(w_k).
\]
Taking norms and using the bounds above yields
\[
    \|w_{k+1}-w^*\|_2
    \le q_0(\ell)\,\|w_k-w^*\|_2
        + \eta\,\bigl\|I - \tau\alpha_k^{(\ell)}\bigr\|_2\,
           \|R(w_k)\|_2.
\]
By Lemma~\ref{lem:lm-bounded-MN} and Corollary~\ref{cor:lm-bounded-alpha}
there exists $C_I(\ell)>0$ such that
\(
    \|I - \tau\alpha_k^{(\ell)}\|_2 \le (1+ \tau C_\alpha (\ell))=:C_I(\ell)
\)
for all $k$, hence
\[
    \|w_{k+1}-w^*\|_2
    \le q_0(\ell)\,\|w_k-w^*\|_2
       + \eta C_I(\ell) C_R \,\|w_k-w^*\|_2^2.
\]

Choose $0<r(\ell)\le r_0$ such that
\[
    \eta C_I(\ell) C_R\, r(\ell)
    \le \frac{1 - q_0(\ell)}{2}.
\]
If $\|w_k-w^*\|_2\le r(\ell)$, the previous inequality gives
\[
    \|w_{k+1}-w^*\|_2
    \le \Bigl(q_0(\ell) + \eta C_I(\ell)C_R r(\ell)\Bigr)\|w_k-w^*\|_2
    \le q(\ell)\,\|w_k-w^*\|_2,
\]
where $q(\ell) := \bigl(1+q_0(\ell)\bigr)/2 \in (0,1)$.
Starting from $w_0\in B_{r(\ell)}(w^*)$ and arguing by induction we obtain
that $w_k\in B_{r(\ell)}(w^*)$ for all $k$ and
\(
    \|w_k-w^*\|_2 \le q(\ell)^k\|w_0-w^*\|_2
\),
which proves both (1) and (2).
\end{proof}

The theorem shows that, for any fixed history length $\ell$, LM--LRSGA
inherits the local linear convergence guarantees of LRSGA, up to a change in
the admissible stepsizes $\tau$ and $\eta$.
The truncation of the memory is reflected only in the constant $\delta^{(\ell)}$ in
Corollary~\ref{cor:lm-bounded-alpha} and hence in the admissible range for $\tau$.

In the next section, we present experiments on generative adversarial networks (GANs), benchmarking against the widely used Adam optimizer, to empirically assess the advantages of the Limited–Memory LRSGA (LM–LRSGA) variant.

\section{Experimental results} \label{sec:experiments}
Using two benchmark datasets, we have measured the performance of our algorithm comparing with the well known Adam optimizer. Regarding the constraint of using a competitive environment, we choose the Generative Adversarial Networks (GANs) model that has two models trained in a competitive manner, the performance of an optimizer is determined by its ability to balance the dynamics between two adversarial objectives: the generator goal of producing realistic outputs and the discriminator goal of correctly distinguishing between real and generated samples. This setting is inherently unstable, as improvements in one model performance can degrade the other, leading to oscillations or divergence if the optimization is not well controlled. In the standard GAN formulation, the discriminator $D$ seeks to maximize the classification accuracy between real data $x \sim p_{\text{data}}$ and generated samples $G(z)$ from latent variables $z \sim p_{z}$ (i.e., randomly sampled from a uniform distribution), while the generator $G$ aims to minimize the same objective to ``fool'' the discriminator. Formally, the discriminator loss is given by

$$
\mathcal{L}_{D} = -\mathbb{E}_{x \sim p_{\text{data}}} \left[ \log D(x) \right] - \mathbb{E}_{z \sim p_{z}} \left[ \log \left( 1 - D(G(z)) \right) \right],
$$

and the generator loss by

$$
\mathcal{L}_{G} = -\mathbb{E}_{z \sim p_{z}} \left[ \log D(G(z)) \right] + \lambda \| \theta_G \| 
$$

where the $\lambda$ term weights the penalty component for the G parameters.

An effective optimizer in this context must maintain a delicate equilibrium, ensuring that neither the generator nor the discriminator overpowers the other, thereby facilitating stable training and yielding high-quality, realistic outputs over time.

We have analyzed the shape of the losses to understand the ability of the optimizer of minimize and maximize the training objectives. Also, we present a quality metric to understand how the trained model is able to generate realistic image. In this context, the term ``realistic'' refers to ability of generating images similar to real data.
The Fréchet Inception Distance (FID) score is a widely used metric for evaluating the quality of images generated by generative models such as GANs. It measures the similarity between the distribution of generated images and that of real images by comparing the mean and covariance of feature representations extracted from a pretrained Inception network. Mathematically, lower FID values indicate that the generated images are more similar to real ones in terms of both visual quality and diversity, with zero representing a perfect match.

\subsection{Experiment on Mnist}
The first experiment is conducted on the MNIST dataset, a classic benchmark in machine learning and computer vision, consisting of 70,000 grayscale images of handwritten digits from 0 to 9. Each image is $28 \times 28$ pixels and centered within the frame to standardize positioning. The data set is divided into 60,000 training samples and 10,000 test samples, and its simplicity, small size, and well-defined labels have made it a fundamental starting point for testing and comparing classification algorithms. We used only the training split to train the GAN.

\subsubsection{Results}
\begin{figure}[htbp]
    \centering
    \begin{subfigure}[b]{0.48\textwidth}
        \centering
        \includegraphics[width=\textwidth]{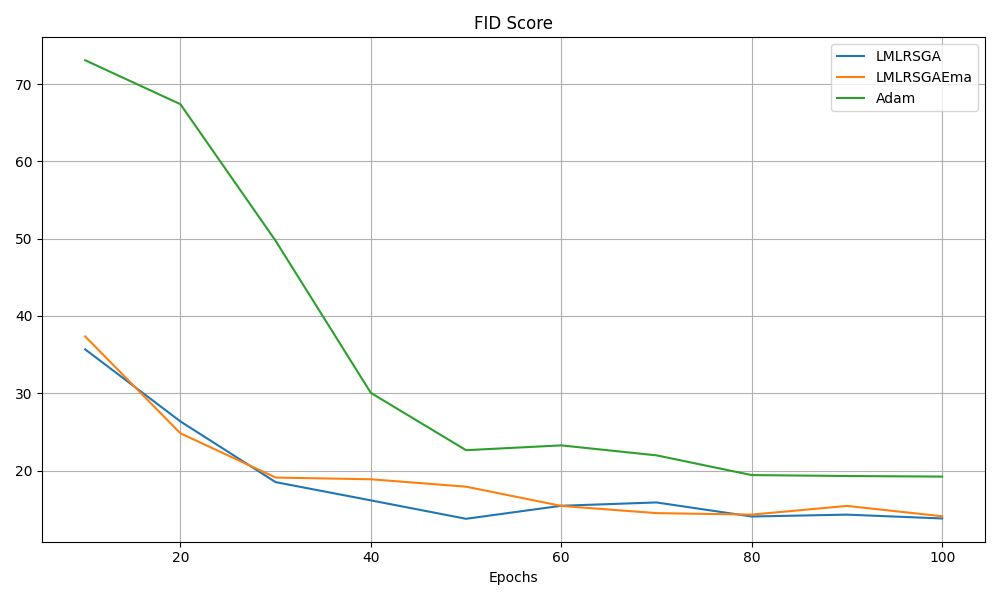}
        \caption{FID scores on MNIST dataset.}
        \label{fig:mnist_fid}
    \end{subfigure}
    \hfill
    \begin{subfigure}[b]{0.48\textwidth}
        \centering
        \includegraphics[width=\textwidth]{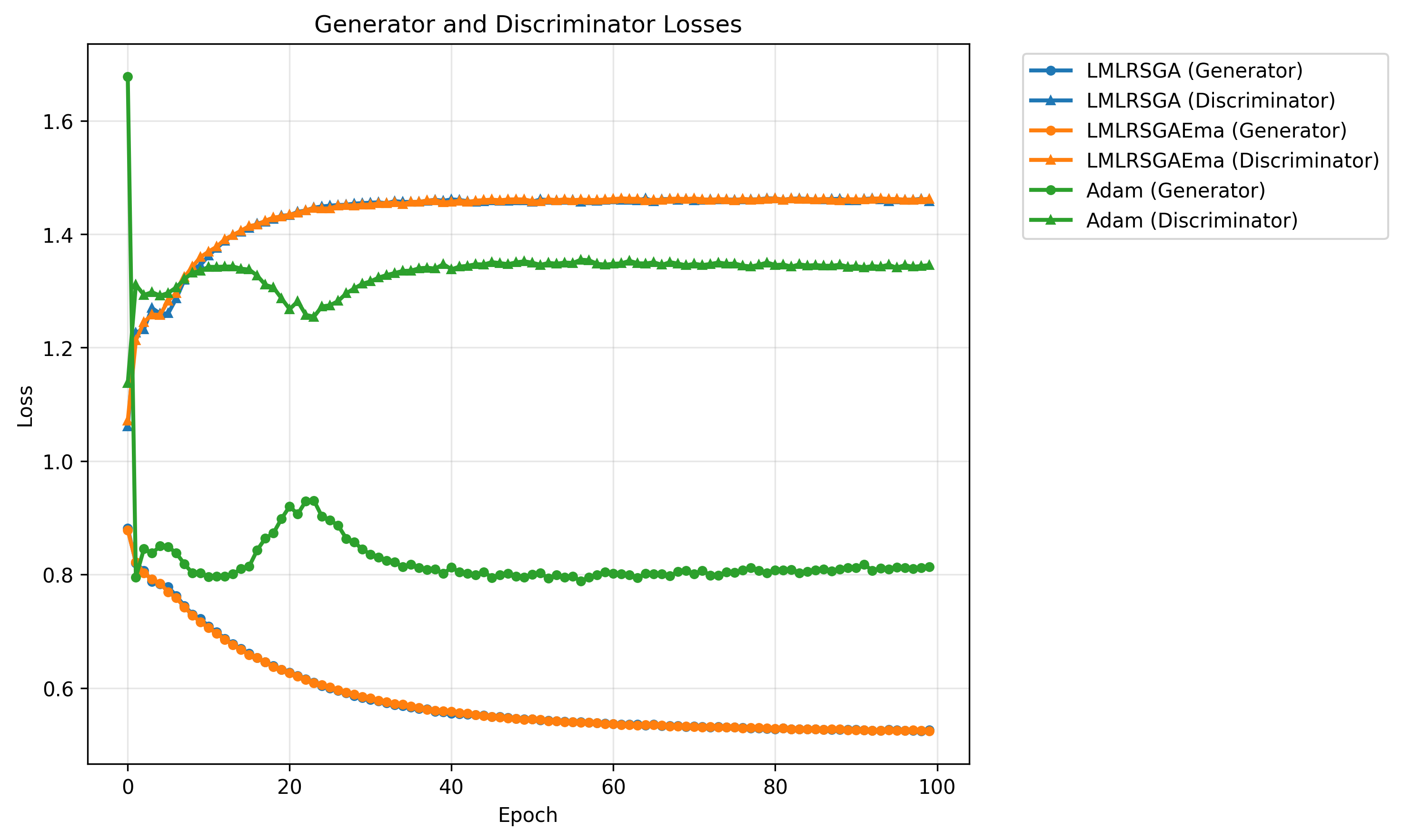}
        \caption{Training losses on MNIST dataset.}
        \label{fig:mnist_loss}
    \end{subfigure}
    
    \caption{Comparison of LM-LRSGA, LM-LRSGAEma and Adam optimizers on MNIST datasets.}
    \label{fig:mnist_comparison}
\end{figure}

Figure \ref{fig:mnist_fid} presents the Fréchet Inception Distance (FID) scores over training epochs on the MNIST dataset. It compares our proposed optimizer LM-LRSGA, and its exponential moving average variant (LRSGAEma) against the widely used Adam optimizer. Both LM-LRSGA and LM-LRSGAEma achieve substantially lower FID scores throughout training than Adam, demonstrating superior sample quality. While Adam begins with a significantly higher FID score (above 70) and gradually decreases to around 19 after 100 epochs, LM-LRSGAand LM-LRSGAEma start with lower scores ($\sim 40$) and converge much faster. Notably, LM-LRSGAEma consistently outperforms LM-LRSGA, yielding smoother convergence and the lowest final FID score of approximately 13. These results highlight the effectiveness of our optimizer in stabilizing generative model training and producing higher-quality samples compared to the conventional Adam baseline.

Figure \ref{fig:mnist_loss} shows the generator and discriminator loss trajectories on the MNIST dataset when trained with LM-LRSGA, LM-LRSGAEma, and Adam. The results indicate that LM-LRSGA and LM-LRSGAEma lead to more stable and balanced adversarial training than Adam. With Adam, the generator loss remains consistently lower ($\sim1.3$) while the discriminator stabilizes around $\sim 0.8$, reflecting a tendency toward imbalance and potential mode collapse. In contrast, LM-LRSGA and LM-LRSGAEma exhibit more synchronized dynamics: the generator loss converges to $\sim1.4$, while the discriminator loss steadily decreases toward $\sim0.5$, suggesting healthier competition between the two networks. Notably, LM-LRSGAEma achieves the smoothest convergence, stabilizing both losses earlier and exhibiting minimal oscillation. These findings align with the FID results, demonstrating that our proposed optimizer fosters more stable GAN training and enhances generative performance relative to the Adam baseline.

\subsubsection{Sensitivity analysis}
In this section, we analyze the role of the main parameters of the LM-LRSGAEma algorithm. We studied how the model's ability to generate realistic images changes when we vary the value of a parameter.

\begin{figure}[htbp]
    \centering
    \begin{subfigure}[b]{0.48\textwidth}
        \centering
        \includegraphics[width=\textwidth]{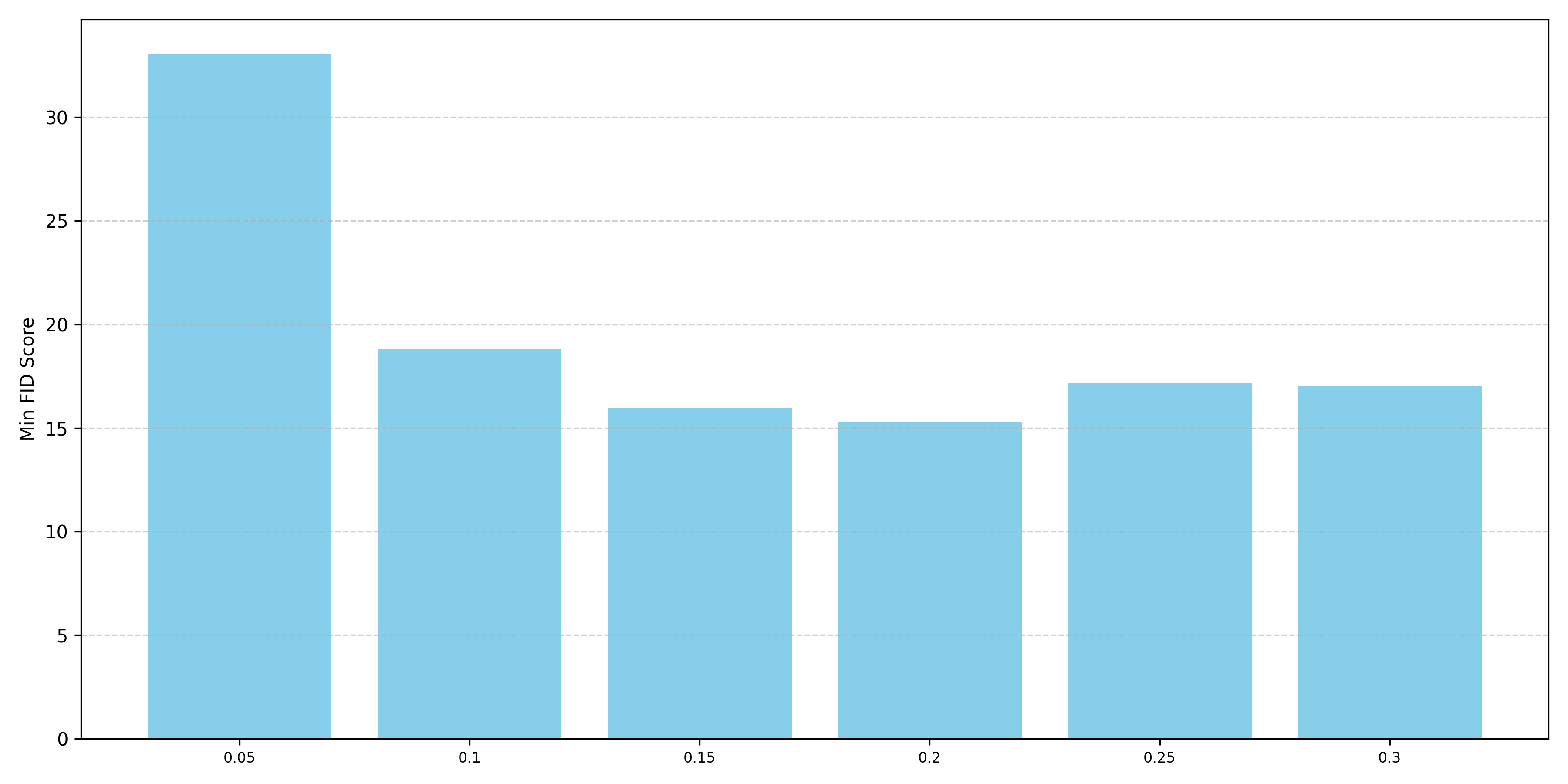}
        \caption{Param $\eta$}
        \label{fig:mnist_eta}
    \end{subfigure}
    \hfill
    \begin{subfigure}[b]{0.48\textwidth}
        \centering
        \includegraphics[width=\textwidth]{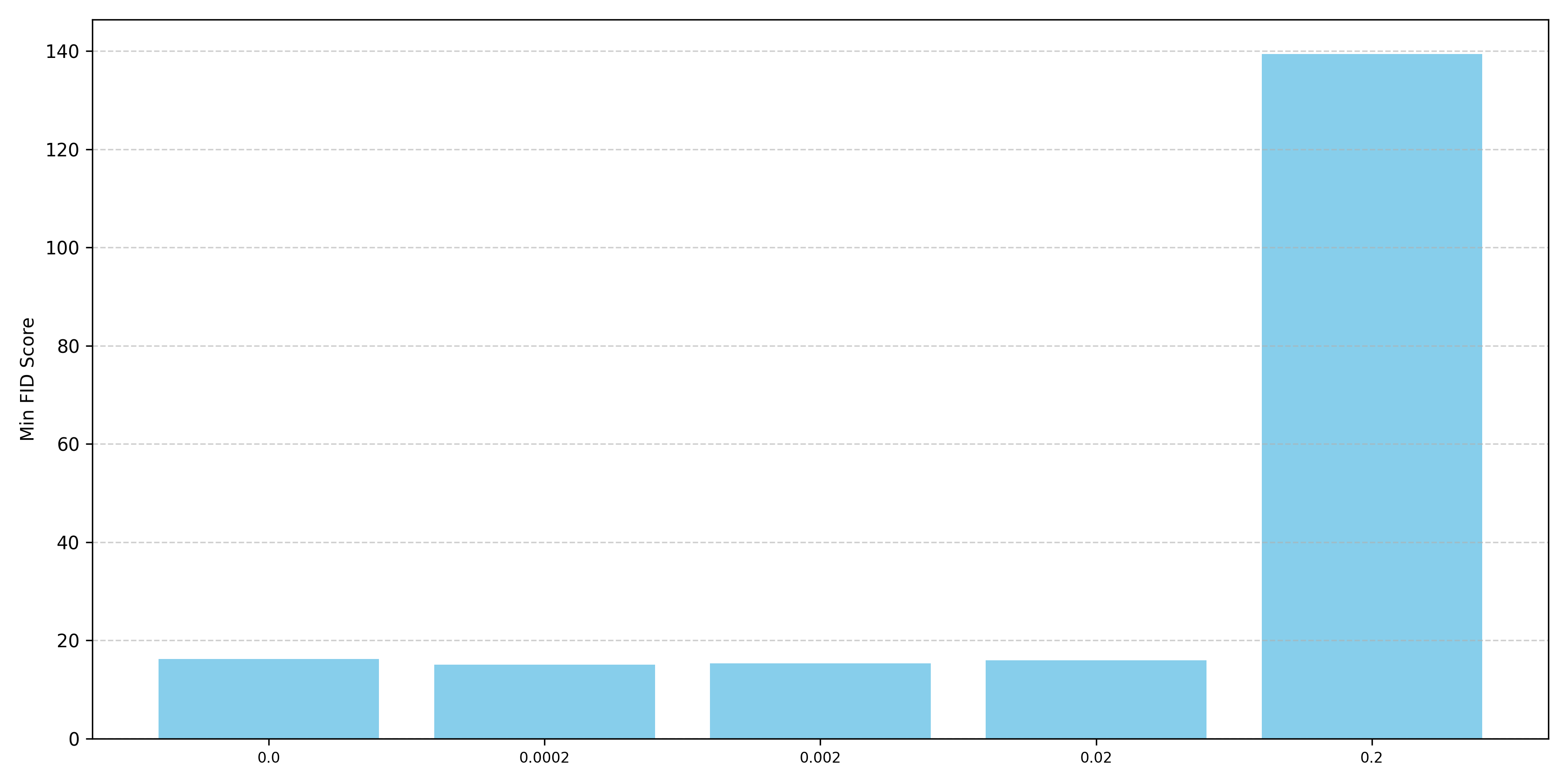}
        \caption{Param $\tau$}
        \label{fig:mnist_tau}
    \end{subfigure}

    \begin{subfigure}[b]{0.48\textwidth}
        \centering
        \includegraphics[width=\textwidth]{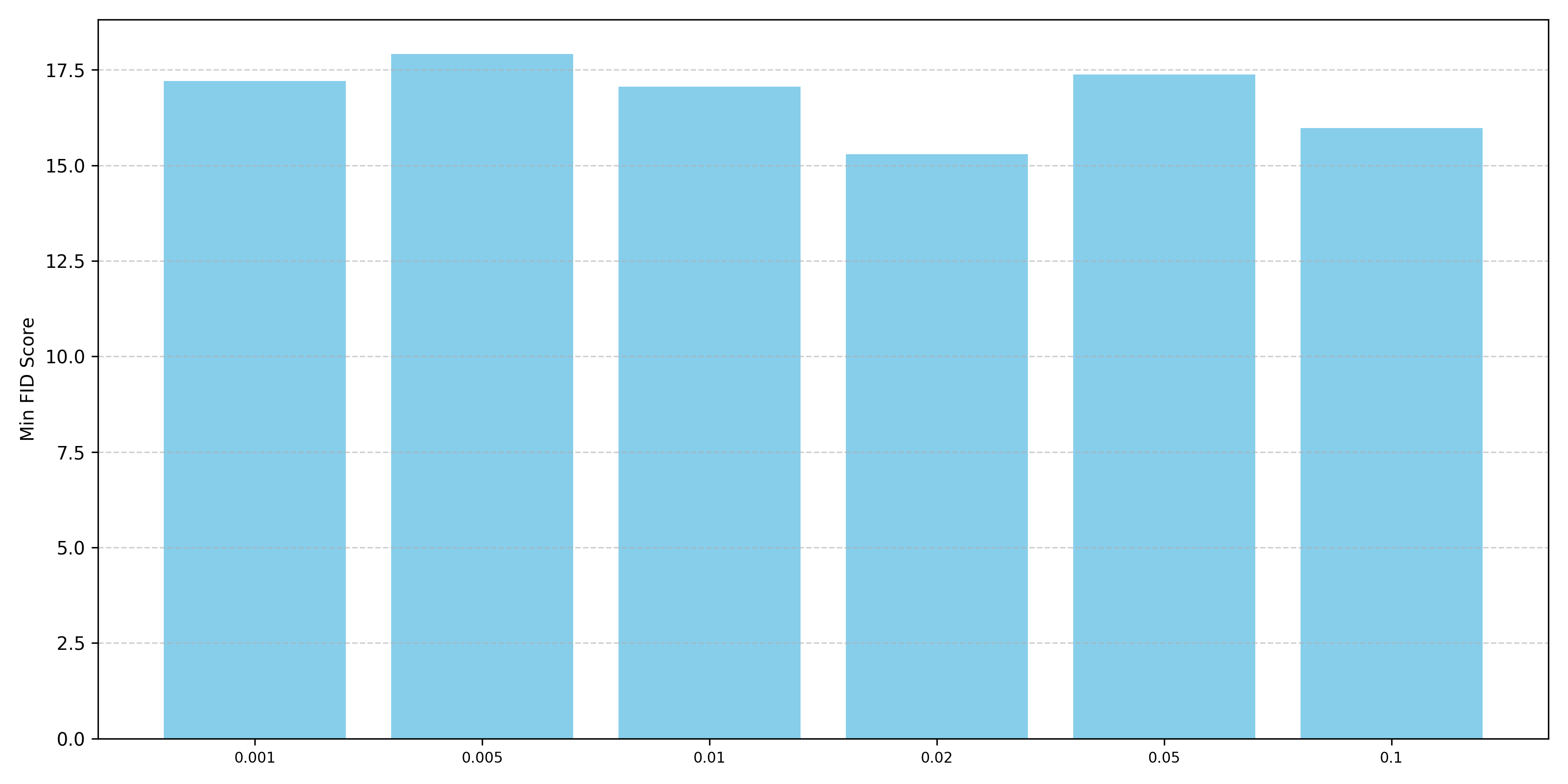}
        \caption{Param $\varepsilon$}
        \label{fig:mnist_eps}
    \end{subfigure}
    \hfill
    \begin{subfigure}[b]{0.48\textwidth}
        \centering
        \includegraphics[width=\textwidth]{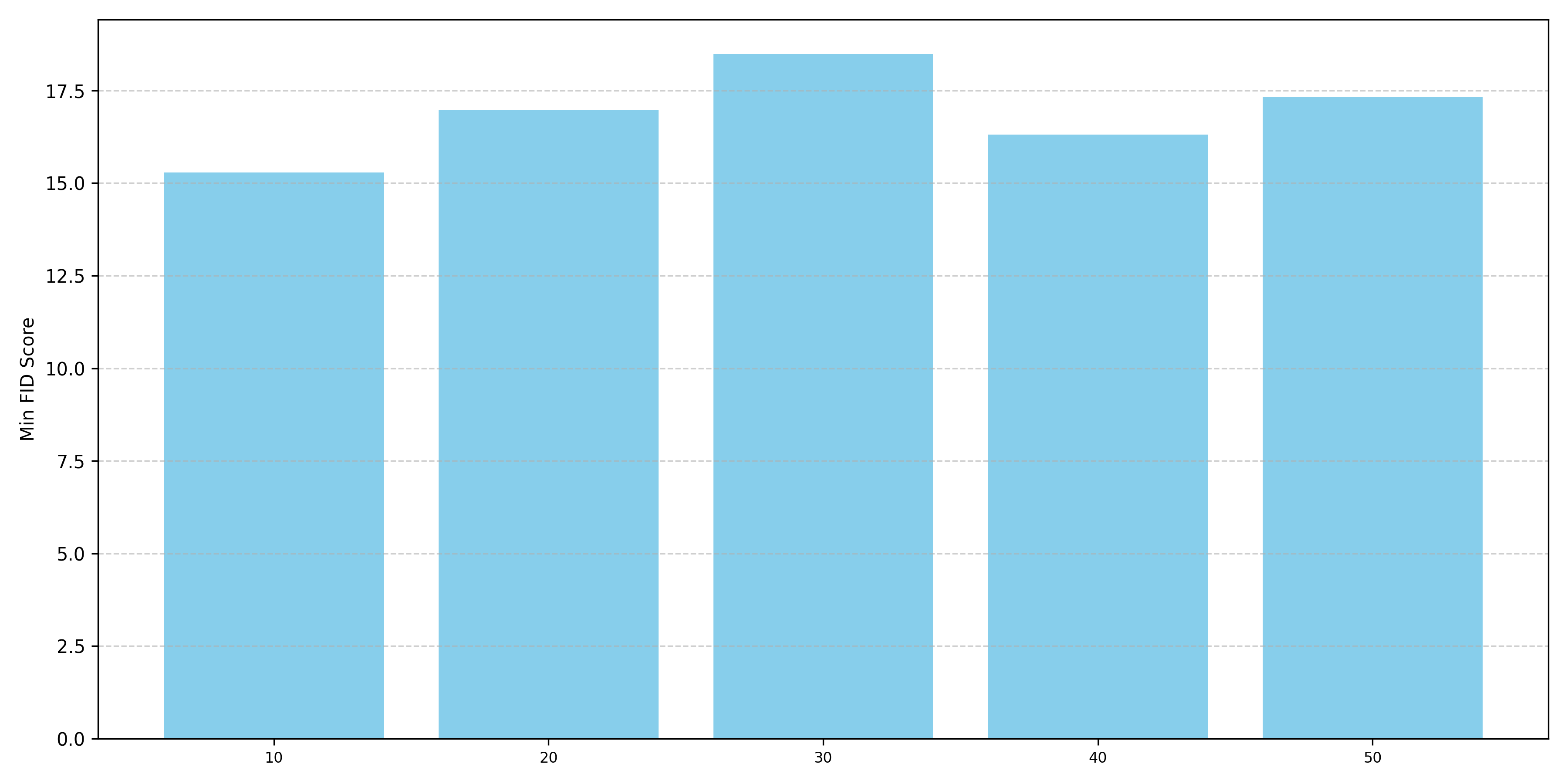}
        \caption{Param Max-History}
        \label{fig:mnist_history}
    \end{subfigure}
    
    \caption{Sensitivity analysis about LM-LRSGAEma parameters.}
    \label{fig:sensitivity_analysis}
\end{figure}

Figure \ref{fig:mnist_eta} presents the hyperparameter sensitivity analysis $\eta$ on the minimum Fréchet Inception Distance (FID) score achieved during the experiments. The bar chart shows results for $\eta$ values ranging from 0.05 to 0.3. Performance varies notably across settings, with the highest (worst) FID score observed at $\eta = 0.05$, indicating degraded generative quality in this regime. The best results are obtained for $\eta$ values between 0.15 and 0.2, yielding minimum FID scores around 15, suggesting an optimal trade-off between learning stability and image fidelity. Larger $\eta$ values beyond 0.25 slightly worsen the FID, though still outperforming the smallest $\eta$ tested. This pattern highlights the importance of tuning $\eta$ for maximizing GAN performance.

Figure \ref{fig:mnist_tau} reports the hyperparameter sensitivity analysis for the parameter $\tau$. The results are shown for $\tau$ values from 0 to 0.2. The 0 value is an experiment in which the component is ignored. Across this range, the minimum FID scores vary only slightly except for value 0.2, indicating that the model performance is relatively insensitive to $\tau$ within the tested limits. The best generative quality is observed at $\tau = 0.002$, both achieving FID scores just above 16, while the highest scores (worst performance) occur around $\tau = 0.2$. In general, the flat trend for values below the 0.02 suggests that $\tau$ has a limited effect on GAN performance in this setup only when it is weighted below a threshold.

 Figures \ref{fig:mnist_eps} and \ref{fig:mnist_history} show the results of the experiments to evaluate the effect of $\varepsilon$ and \textit{max\_history} parameters. The charts show only negligible differences between the settings. For the $\varepsilon$ parameter, we expect its impact to be minimal, as it is effectively offset by the residuals observed between consecutive optimization iterations. As for the \textit{max\_history} parameter, we consider a value of 10 sufficient for this dataset to accurately estimate the correction to be applied and the resource usage required to execute the experiment.

\subsection{Experiment on Fashion Mnist}
In this section we analyze the results about the Fashion-MNIST dataset, it is a widely used benchmark for image classification tasks, designed as a more challenging drop-in replacement for the original MNIST handwritten digits. It consists of 70,000 grayscale images of size $28 \times 28$ pixels, split into 60,000 training samples and 10,000 test samples. We used only the training split to train the GAN models. Each image depicts a single fashion item from one of 10 categories, including T-shirts, trousers, pullovers, dresses, etc. The dataset's compact size and diversity make it a standard testbed for evaluating machine learning algorithms in vision research.

\begin{figure}[htbp]
    \centering
    \begin{subfigure}[b]{0.48\textwidth}
        \centering
        \includegraphics[width=\textwidth]{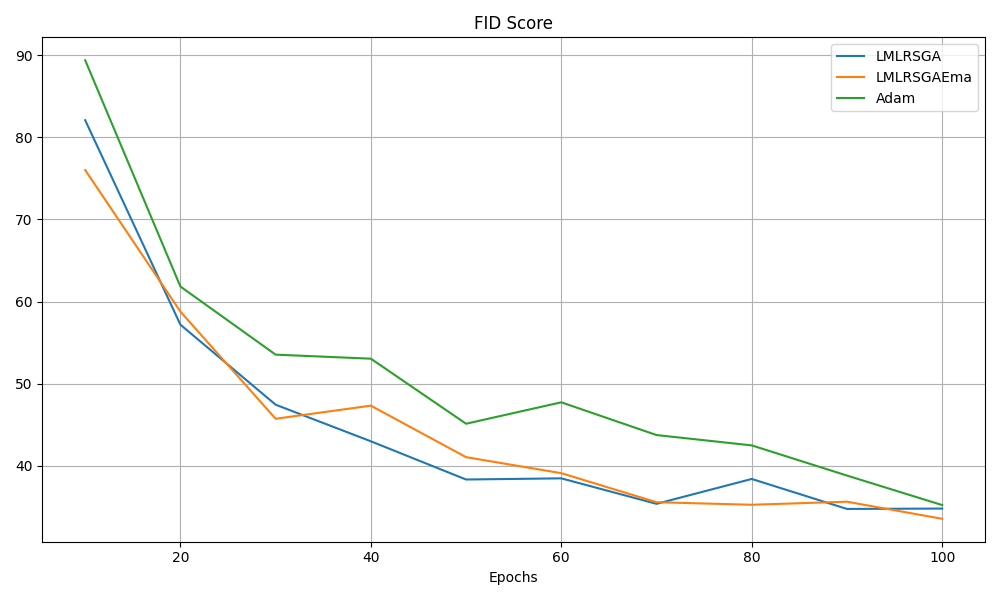}
        \caption{FID scores on Fashion dataset.}
        \label{fig:fashion_fid}
    \end{subfigure}
    \hfill
    \begin{subfigure}[b]{0.48\textwidth}
        \centering
        \includegraphics[width=\textwidth]{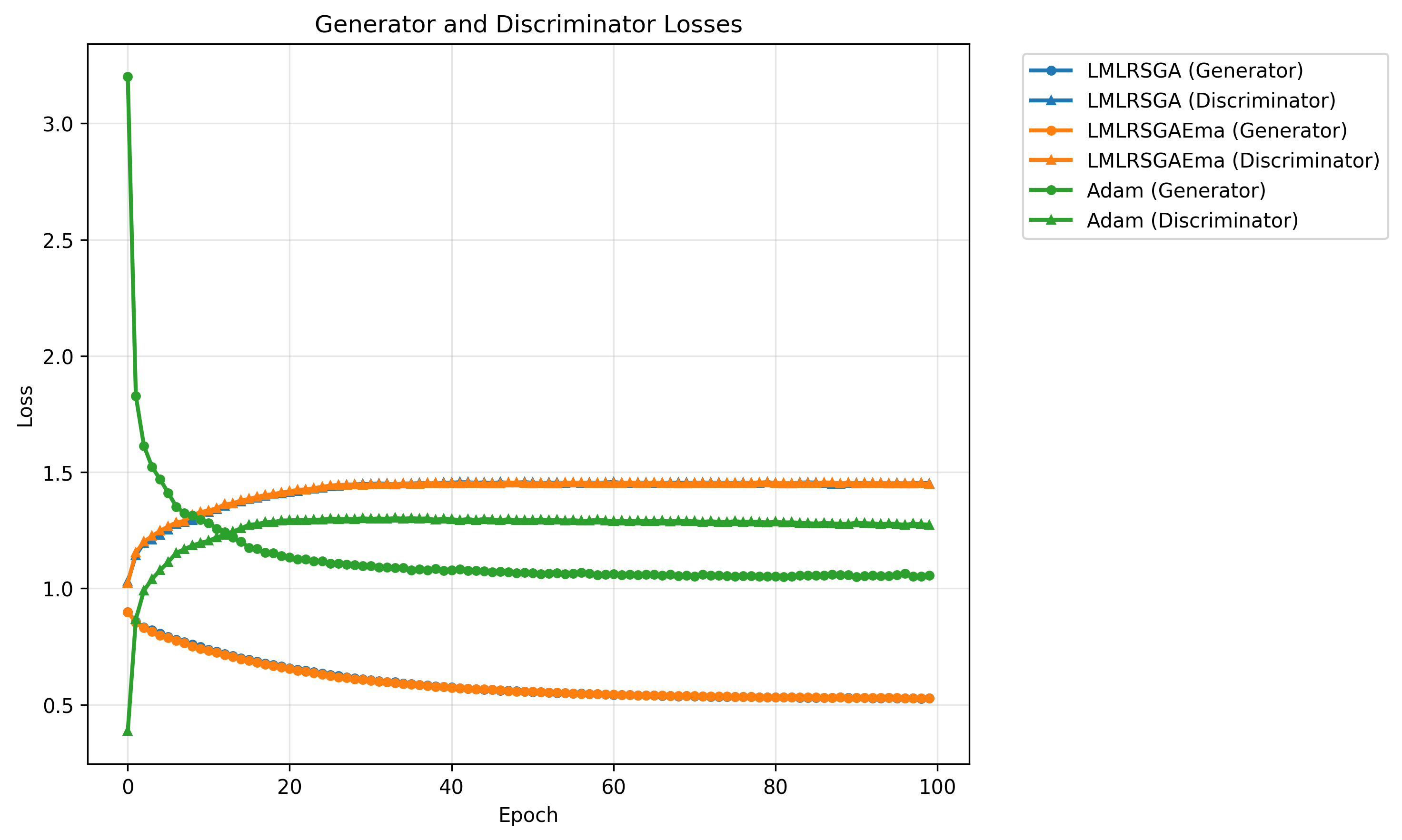}
        \caption{Training losses on Fashion dataset.}
        \label{fig:fashion_loss}
    \end{subfigure}
    
    \caption{Comparison of LM-LRSGA, LM-LRSGAEma, and Adam optimizers on Fashion datasets.}
    \label{fig:fashion_comparison}
\end{figure}

Figure \ref{fig:fashion_fid} shows the evolution of the FID scores over 100 training epochs for the proposed LM-LRSGA and LM-LRSGAEma optimizer and the Adam baseline optimizer on the Fashion dataset. The FID score decreases sharply for all methods during the initial epochs, dropping from 90 to approximately 55–60 by epoch 20.
LM-LRSGA and LM-LRSGAEma consistently achieves lower FID scores for all epochs. Adam achieves comparable performance on the last epoch. We have limited the analysis to 100 epochs to maintain coherence with Mnist experiments.

Figure \ref{fig:fashion_loss} shows the training losses of Adam, LM-LRSGA and LM-LRSGAEma. The same behavior of Mnist dataset is observed in this experiment. The two losses converge to a stable value around the epoch $\sim 30$ for LM-LRSGA and LM-LRSGAEma, while Adam converges to some epochs before. The distance between the two losses confirms the ability of LM-LRSGA and LM-LRSGAEma to better deal than Adam.  

\section{Spectral and Stability Analysis}
\label{sec:spectral_analysis}

To complement the quantitative evaluation based on FID, we provide an empirical characterization of the dynamical stability of the adversarial training process.
We first describe the spectral framework and the quantitative metrics adopted for the analysis, then present the experimental setup and the main results.

\subsection{Spectral and Stability Metrics}
\label{subsec:spectral_metrics}

The spectral and stability analysis relies on an empirical estimation of the local Jacobian of the joint generator--discriminator dynamics, following the theoretical framework for differentiable games introduced in \cite{mescheder2018gan,balduzzi2018mechanics,fiez2021spectral}. 
The procedure combines dimensionality reduction, spectral decomposition, and power--spectral analysis to characterize the linearized dynamics near equilibrium.

\paragraph{Jacobian spectrum estimation.}
Given the temporal sequence of parameter vectors 
$w_k = [\theta_G^{(k)}, \theta_D^{(k)}]^\top \in \mathbb{R}^{m+n}$, 
the discrete dynamics are approximated locally by a linear operator
$w_{k+1} \approx A w_k$, 
where $A$ is obtained by solving a least-squares regression between successive iterates. 
To improve conditioning and isolate the dominant modes of variation, the trajectories are projected onto a reduced subspace using truncated Singular Value Decomposition (SVD). 
This projection corresponds to the standard low-rank approximation 
$X \approx U_r \Sigma_r V_r^\top$, 
which filters numerical noise and captures the principal components of the joint update. 
The reduced operator $A^\star = (Y V_r)\Sigma_r^{-1}U_r^\top$ is then analyzed through an eigen decomposition performed via a block-Lanczos iterative scheme, which allows an efficient approximation of the leading eigenvalues without explicitly forming dense matrices. 
The eigenvalues $\{\lambda_i\}$ of $A^\star$ provide an empirical estimate of the local Jacobian spectrum, and the \textit{spectral radius} is defined as $\rho = \max_i |\lambda_i|$.

\paragraph{Stability regime classification.}
The notions of \textit{stable}, \textit{marginally stable}, and \textit{unstable} dynamics follow the standard classification adopted in the spectral analysis of differentiable and adversarial games \cite{mescheder2018gan,balduzzi2018mechanics,fiez2021spectral}. 
In discrete-time settings, linear stability is determined by the spectral radius of the estimated Jacobian:
$\rho < 1$ indicates convergence (stable regime), 
$\rho = 1$ corresponds to marginal stability, 
and $\rho > 1$ implies divergence. 
To account for numerical noise and finite-length trajectories, a tolerance band of $\epsilon = 0.15$ is applied. 
Cases with $\rho \in [1-\epsilon,\, 1+\epsilon]$ are further analyzed by inspecting the frequency content of the loss trajectories. 
This inspection is based on the Fourier transform of the temporal loss signals, where the tolerance $\epsilon$ effectively separates slowly varying components from high-frequency oscillatory modes. 
The spectral energy of the oscillatory band is estimated using \textit{Welch’s method} \cite{welch1967psd}, which averages modified periodograms over overlapping windows to reduce the variance of the power spectral density (PSD) estimate. 
In practice, we compute a high-frequency power index for both the generator and the discriminator; when these indices are of order $10^{0}$, the dynamics are classified as \textit{marginally stable}, whereas regimes with significantly larger high-frequency power (i.e., well above $10^{0}$) are deemed \textit{unstable}.
This refinement is consistent with the interpretation of lightly damped rotational dynamics commonly observed in adversarial training systems \cite{fiez2021spectral}.

\paragraph{Auxiliary stability metrics.}
In addition to the spectral radius, several scalar stability indicators are evaluated to summarize the temporal smoothness and robustness of the optimization process. 
The \textit{loss stability} index measures the local variability of the generator loss through the inverse of its standard deviation across updates, providing a smoothness score between 0 and 1. 
The \textit{mode-collapse trend} quantifies the time derivative of the variance of the generator parameters and serves as a proxy for distributional diversity, following the methodology proposed in \cite{lucic2018ganmetrics}. 
Finally, the \textit{global stability index} computes the inverse standard deviation of parameter trajectories within sliding windows, functioning as an empirical analogue of a local Lyapunov factor. 
All metrics are normalized in $[0,1]$, where larger values indicate smoother and more stable dynamics.

\subsection{Experimental Setup}
\label{sec:exp_setup}

All experiments were conducted on the same MNIST GAN architecture used for the FID evaluation, ensuring full comparability of results. 
The training dynamics of the three optimizers (Adam, LM-LRSGA, and LM-LRSGAEma) were analyzed under two \textit{fair} learning-rate regimes:
\begin{enumerate}[label=(\roman*)]
    \item \textbf{LM-LRSGA-optimal setup}: $\eta = 0.2$, $\tau = 0.002$, corresponding to the FID best-performing configuration for both LM-LRSGA and LM-LRSGAEma;
    \item \textbf{Adam-optimal setup}: $\eta = 0.001$, $\tau = 0.002$, corresponding to the learning rate that maximized the FID performance for Adam.
\end{enumerate}

For each optimizer, the joint trajectories of generator and discriminator parameters were recorded over the entire training process. 
Parameter updates were collected at uniform intervals and concatenated to form the sequence 
$w_k = [\theta_G^{(k)}, \theta_D^{(k)}]^\top$. 
To reduce dimensionality and numerical noise, the trajectories were projected onto a 40-dimensional subspace obtained by truncated SVD. 
The empirical Jacobian operator was then estimated from these projected trajectories using the least-squares formulation described in Section~\ref{subsec:spectral_metrics}. 
The dominant eigenvalues of the resulting operator were computed through a block-Lanczos iterative decomposition, which efficiently approximates the leading part of the spectrum without constructing dense matrices.

To assess oscillatory behavior, the temporal loss signals of the generator and discriminator were further analyzed in the frequency domain. 
The power spectral density (PSD) was estimated using \textit{Welch’s method} \cite{welch1967psd}, as implemented in the \texttt{scipy.signal} module. 
High-frequency power was integrated beyond a cutoff threshold to quantify the spectral energy of oscillations. 
This quantity was used as a secondary indicator to distinguish marginally stable from unstable regimes in cases where $\rho \approx 1$.

A tolerance parameter of $\epsilon = 0.15$ was adopted in all experiments for the stability classification, as discussed in Section~\ref{subsec:spectral_metrics}.
All analyses were performed using Python 3.10 with the \texttt{NumPy}, \texttt{SciPy}, and \texttt{PyTorch} libraries.

\subsection{Results and Discussion}
\label{sec:results_discussion}

\begin{figure*}[t!]
\centering
\begin{subfigure}[t]{0.95\textwidth}
    \centering
    \includegraphics[width=\textwidth]{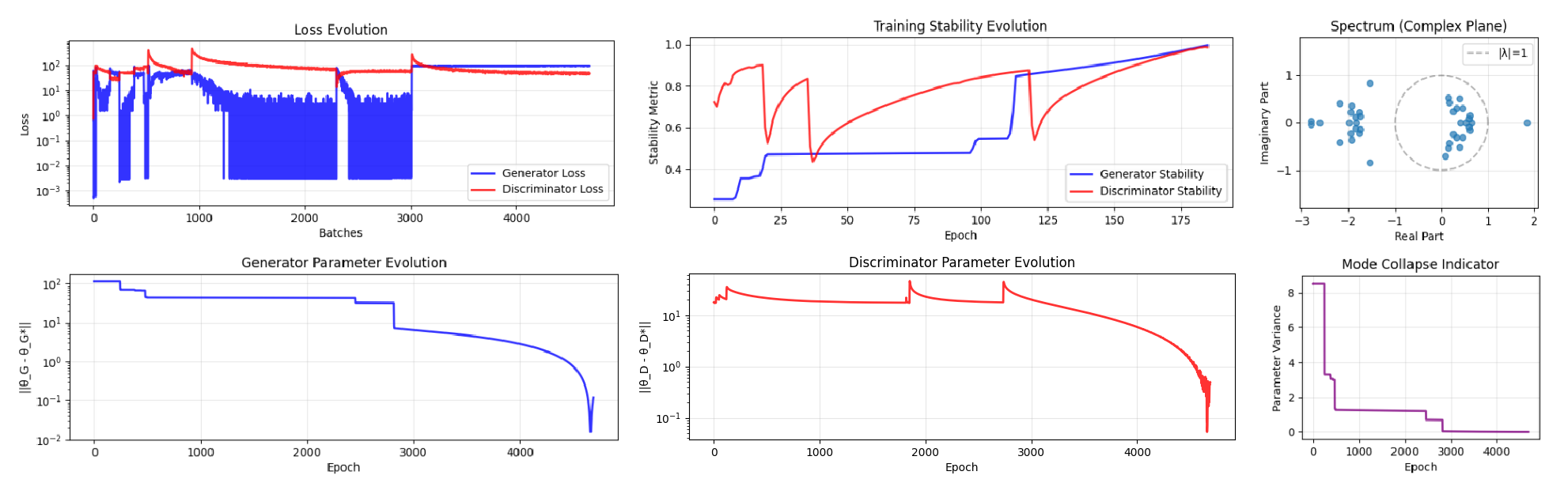}
    \caption{Adam}
    \label{fig:spectral_adam_015}
\end{subfigure}

\vspace{1em}

\begin{subfigure}[t]{0.95\textwidth}
    \centering
    \includegraphics[width=\textwidth]{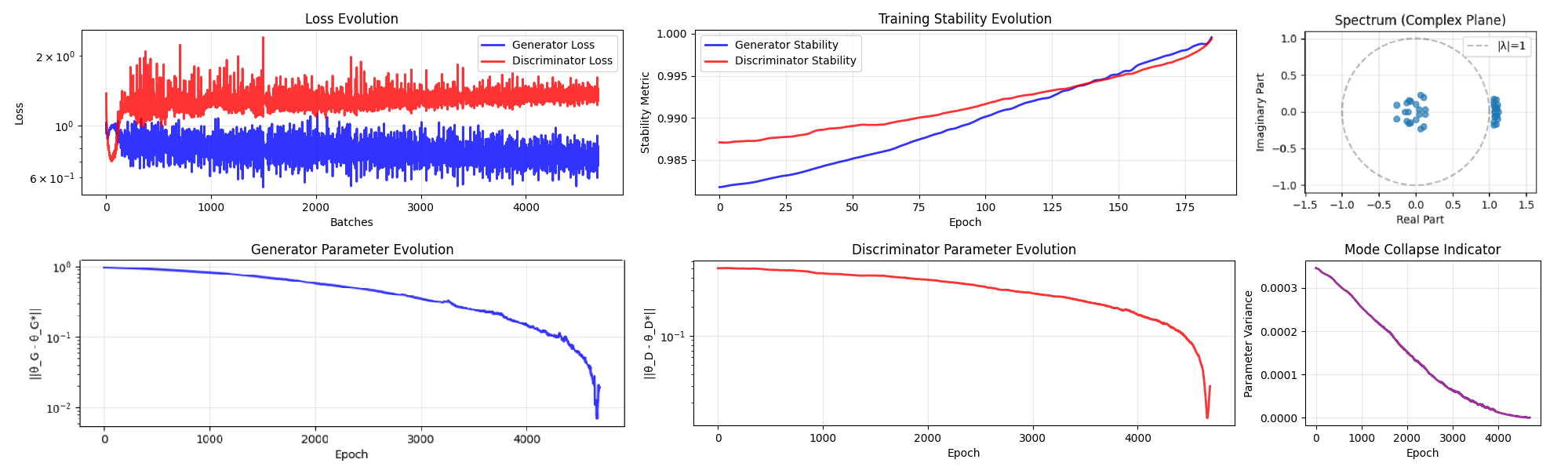}
    \caption{LM-LRSGA}
    \label{fig:spectral_lmlrsga_015}
\end{subfigure}

\vspace{1em}

\begin{subfigure}[t]{0.95\textwidth}
    \centering
    \includegraphics[width=\textwidth]{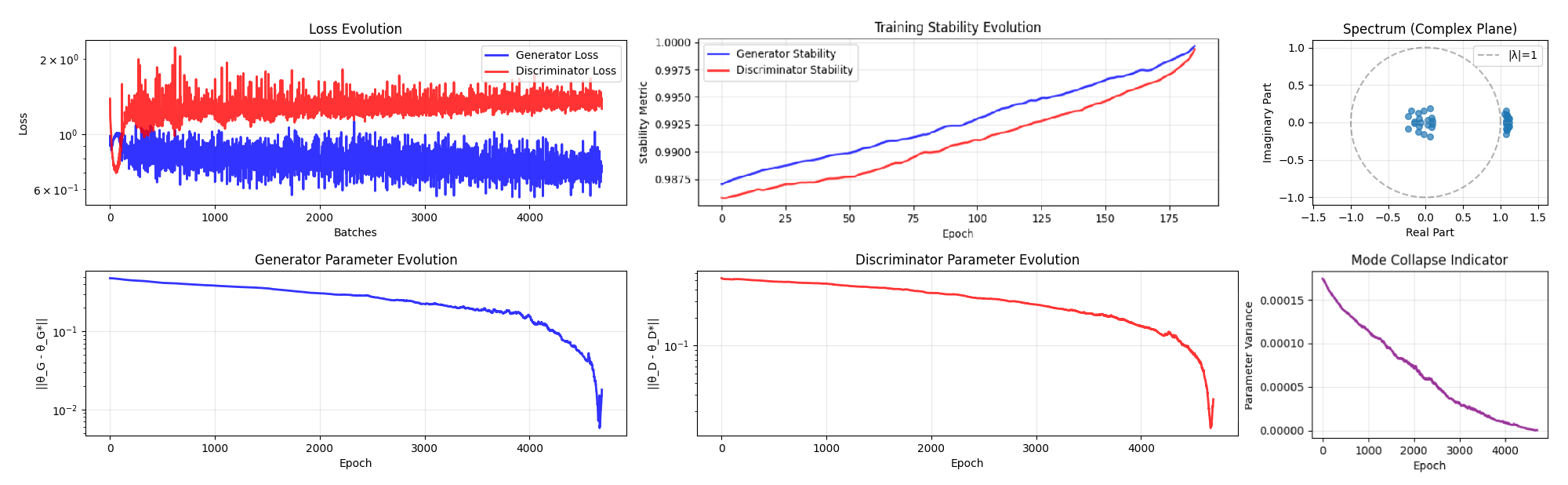}
    \caption{LM-LRSGAEma}
    \label{fig:spectral_lmlrsgaema_015}
\end{subfigure}

\caption{
Spectral and stability diagnostics for Adam, LM-LRSGA, and LM-LRSGAEma with learning rate $\eta = 0.2$. 
Each panel shows (top) loss evolution in logarithmic scale, the training stability metric, and the estimated Jacobian spectrum in the complex plane; 
and (bottom) the generator and discriminator parameter evolution together with the mode-collapse indicator.
}
\label{fig:spectral_results_015}
\end{figure*}

\begin{figure*}[t!]
\centering
\begin{subfigure}[t]{0.95\textwidth}
    \centering
    \includegraphics[width=\textwidth]{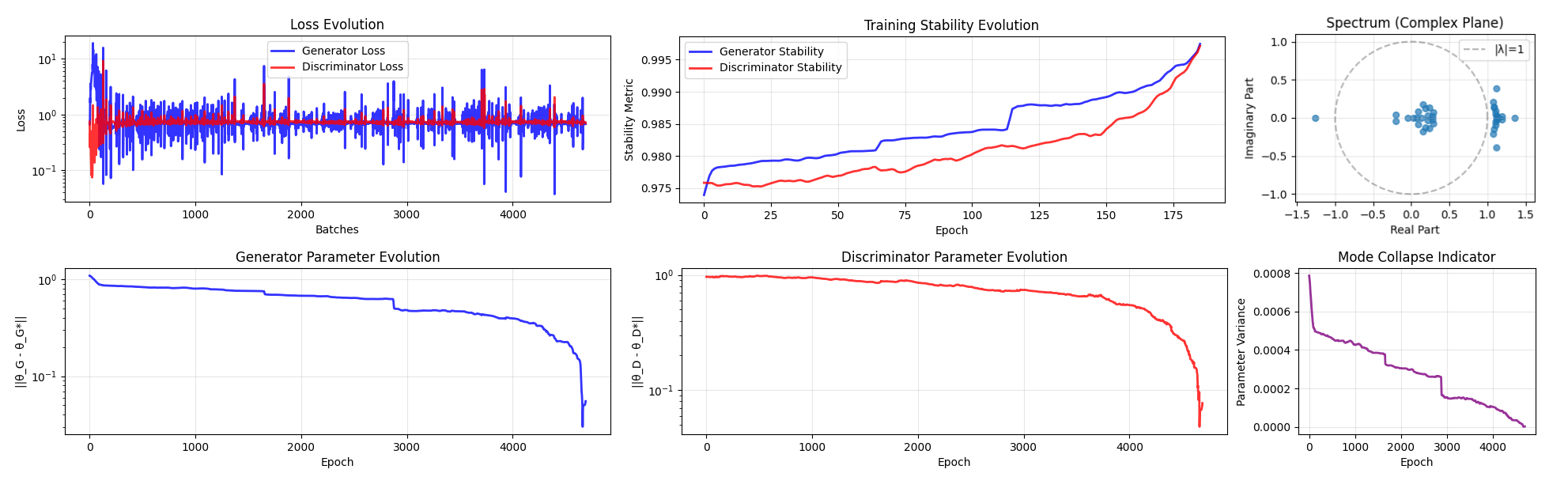}
    \caption{Adam}
    \label{fig:spectral_adam_0001}
\end{subfigure}

\vspace{1em}

\begin{subfigure}[t]{0.95\textwidth}
    \centering
    \includegraphics[width=\textwidth]{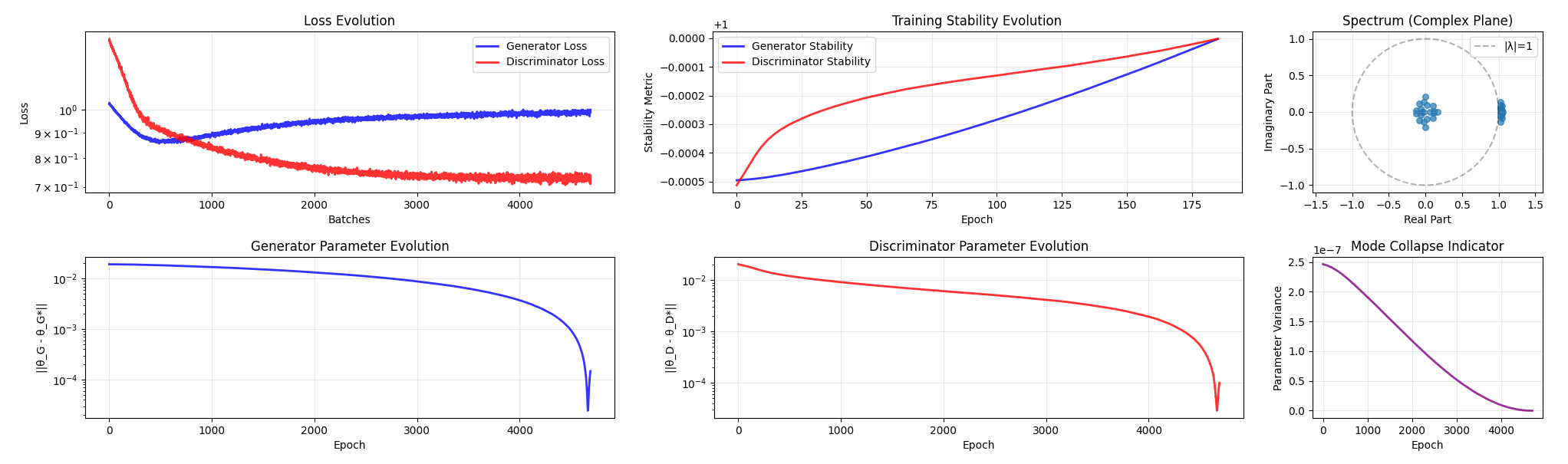}
    \caption{LM-LRSGA}
    \label{fig:spectral_lmlrsga_0001}
\end{subfigure}

\vspace{1em}

\begin{subfigure}[t]{0.95\textwidth}
    \centering
    \includegraphics[width=\textwidth]{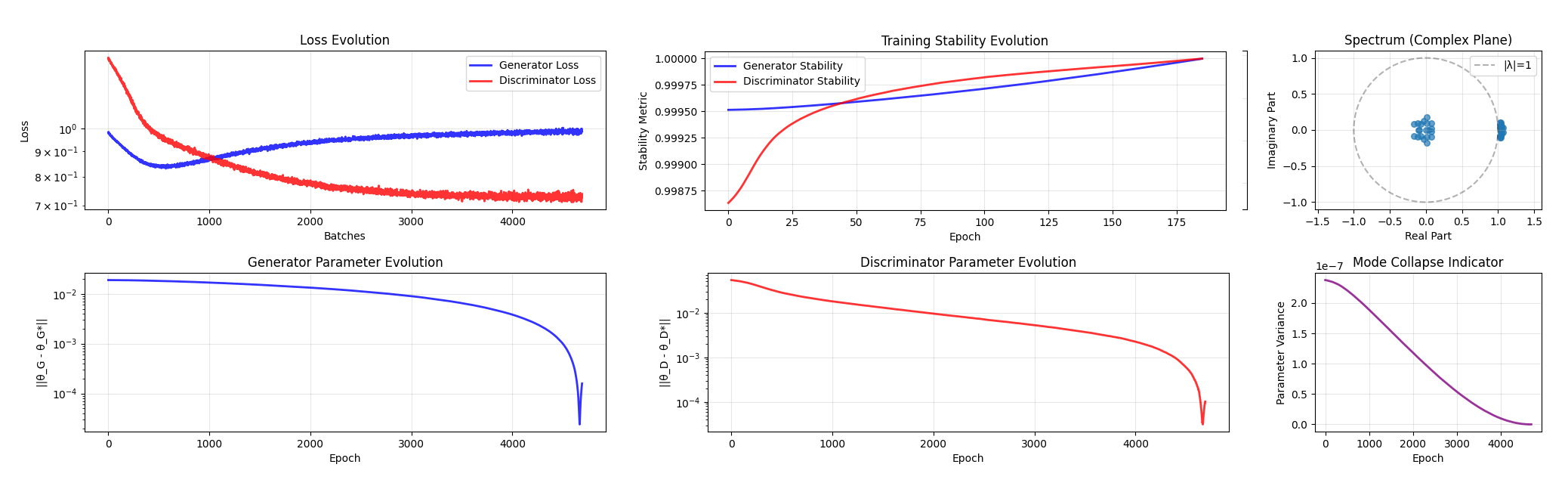}
    \caption{LM-LRSGAEma}
    \label{fig:spectral_lmlrsgaema_0001}
\end{subfigure}

\caption{
Spectral and stability diagnostics for Adam, LM-LRSGA, and LM-LRSGAEma with learning rate $\eta = 0.001$. 
Each panel shows (top) loss evolution in logarithmic scale, the training stability metric, and the estimated Jacobian spectrum in the complex plane; 
and (bottom) the generator and discriminator parameter evolution together with the mode-collapse indicator.
}
\label{fig:spectral_results_0001}
\end{figure*}

Figures~\ref{fig:spectral_results_015} and \ref{fig:spectral_results_0001} provide a comprehensive visualization of the spectral and stability diagnostics obtained for the analyzed optimizers.
The first row reports the evolution of the generator and discriminator losses shown on a logarithmic scale (left), the corresponding training stability metric (center), and the estimated spectrum of the local Jacobian in the complex plane (right). 
The spectral plot displays the eigenvalues computed from the block-Lanczos approximation of the empirical Jacobian, with the dashed unit circle indicating the boundary between contractive and divergent dynamics.
The second row shows the evolution of the generator parameters (left) and discriminator parameters (center), represented by the Euclidean distance from their initialization, and the \textit{mode-collapse indicator} (right).
A nearly monotonic decrease of these distances together with a slow, smooth decay of the variance indicates stable and non-collapsing behavior.
Overall, the combined diagnostics reveal that both LM-LRSGA variants maintain bounded losses, regular parameter trajectories, and eigenvalues concentrated within or near the unit circle, in contrast to the divergent and oscillatory dynamics observed with Adam.
In particular, when using the LM–LRSGA–optimal configuration ($\eta=0.2$, $\tau=0.002$), \textbf{Adam} exhibits a spectral radius of $\rho = 3.06$, 
a value well above the stability band $[1-\varepsilon,1+\varepsilon] = [0.85,1.15]$ (with $\varepsilon = 0.15$), and is therefore classified as \emph{unstable} by the spectral criterion alone. 
For completeness, we also report the high–frequency (HF) power of the loss signals estimated via Welch’s method: both generator and discriminator HF powers are approximately of order $10^{3}$, confirming the presence of strong oscillatory energy.

In contrast, both \textbf{LM–LRSGA} and \textbf{LM–LRSGA–EMA} fall inside the tolerance band, with $\rho = 1.1248$ and $\rho = 1.1335$, respectively. 
For these borderline cases, we apply the Welch refinement: the orders of the HF power for generator and discriminator are $10^{-1}$ and $10^{0}$, respectively, for both LM–LRSGA and LM–LRSGA–EMA—values that are over three orders of magnitude smaller than for Adam and support their classification as \emph{marginally stable}. Consistently, both limited–memory variants maintain nearly constant global stability indices ($\approx 0.99$) and parameter oscillations of order $\mathcal{O}(10^{-4})$, in line with their superior FID scores.

Under the Adam–optimal learning-rate configuration ($\eta=0.001$, $\tau=0.002$),
\textbf{LM–LRSGA} and \textbf{LM–LRSGA–EMA} fall inside the tolerance band with
$\rho \approx 1.06$ and $\rho \approx 1.06$, respectively, and are thus disambiguated via Welch’s PSD.
Their high–frequency (HF) powers are of order $10^{-3}$ for both generator and discriminator, 
for both LM–LRSGA and LM–LRSGA–EMA, which are several orders of magnitude below those observed for Adam.
Accordingly, both limited–memory variants are classified as \emph{marginally stable},
with near–unit global stability indices and parameter oscillations essentially at $\mathcal{O}(10^{-4})$, 
indicating smooth, non–collapsing dynamics.
In contrast, \textbf{Adam} attains $\rho \approx 1.48$ (above $1+\varepsilon = 1.15$), hence is \emph{unstable} by the spectral criterion alone; 
for completeness, its HF powers are of order $10^{-1}$ for the discriminator and $10^{1}$ for the generator, 
confirming residual oscillatory energy in the losses.

\begin{table}[t]
\centering
\caption{Spectral and stability metrics for Adam, LM-LRSGA, and LM-LRSGAEma under two learning-rate regimes. 
All quantities are computed from block-Lanczos spectral estimates and the auxiliary stability indices described in Section~\ref{subsec:spectral_metrics}.}
\label{tab:spectral_results}
\resizebox{0.95\textwidth}{!}{%
\begin{tabular}{llccccc}
\toprule
\textbf{Setting} & \textbf{Optimizer} & $\rho$ &\textbf{HF Power (G/D)} & \textbf{Pred. Stab}
& \textbf{Stab.\ (G/D)} & \textbf{Mode-Collapse} \\
\midrule
%\multirow{3}{*}{$\eta=0.2$, $\tau=0.002$}
\multirow{3}{*}{$\begin{array}{c}
\eta = 0.2 \\
\tau = 0.002
\end{array}$}

 & Adam & 3.06 & $1.1 \times 10^{3} / 1.0 \times 10^{3}$ &
 Unstable & 0.63 / 0.78 & High \\
 & LM-LRSGA & 1.12 & $7.1 \times 10^{-1} / 1.1 \times 10^{0}$ & Marginally Stable & 0.99 / 0.99 & Low \\
 & LM-LRSGAEma & 1.13 & $7.2 \times 10^{-1} / 1.1 \times 10^{0}$ & Marginally Stable  & 0.99 / 0.99 & Low \\
\midrule
%\multirow{3}{*}{$\eta=0.001$, $\tau=0.002$}
\multirow{3}{*}{$\begin{array}{c}
\eta = 0.001 \\
\tau = 0.002
\end{array}$}

 & Adam & 1.48 & $2.1 \times 10^{1} / 6.4 \times 10^{-1}$ &  
 Unstable & 0.99 / 0.98 & Low \\
 & LM-LRSGA & 1.06 & $2.1 \times 10^{-3} / 3.5 \times 10^{-3}$ & Marginally Stable  & 1.00 / 1.00 & Low \\
 & LM-LRSGAEma & 1.06 & $1.7 \times 10^{-3} / 3.2 \times 10^{-3}$ & Marginally Stable & 1.00 / 1.00 & Low \\
\bottomrule

\end{tabular}%
}
\end{table}

Across both learning-rate regimes, the limited-memory variants substantially reduce the spectral radius compared with Adam and maintain marginally stable trajectories even for large step sizes. 
The reduction of $\rho$ and the absence of high-frequency energy in the loss spectra confirm that the LM-LRSGA correction effectively damps the antisymmetric component of the gradient field, yielding dynamics that are nearly contractive in the sense of the spectral criterion. 
These findings provide a consistent spectral explanation for the smoother convergence and improved sample quality obtained by LM-LRSGA in the FID experiments.

\section{Conclusions}
\label{sec:conclusion}

This manuscript extends our prior LRSGA work \cite{lrsga} in four directions: (i) a new per-iteration spectral stability condition for LRSGA near Nash equilibria (Section~3); (ii) an $\varepsilon_x,\varepsilon_y$-regularized mixed-block formulation that avoids forming the full Jacobian updates (Section~4); (iii) a limited-memory variant (LM-LRSGA) based on adapted recursions together with a local linear convergence analysis (Sections~4--5); (iv) an empirical study and spectral diagnostics in adversarial training (Sections~6--7). 
%The overlap with [2025] is limited to standard preliminaries and the definition of LRSGA.
Our LM-LRSGA, building on symplectic gradient corrections and quasi-Newton approximations, retains the stabilizing effect of second-order information while reducing memory and compute via two-loop recursions algorithms. By storing only a short history of curvature pairs, it attains near–second-order stability with first-order–like efficiency.

Our analyses show that LM-LRSGA mitigates the rotational dynamics that typically hinder convergence in adversarial settings. On GAN training with MNIST and Fashion-MNIST, it yields smoother, synchronized loss trajectories, faster convergence, and FID improvement over Adam. These results indicate that limited-memory curvature corrections offer a practical balance of stability, efficiency, and scalability for competitive optimization.

Future work will focus on broadening the applicability of LM-LRSGA beyond two-player games. In particular, we plan to extend the method to competitive (or more generally non-cooperative) settings involving more than two neural architectures, and to assess its empirical and theoretical behavior in these multi-player training scenarios.

LM-LRSGA also supports sustainable, energy-aware machine learning. By avoiding explicit formation and storage of high-dimensional second-order tensors, it reduces memory transfers and expensive operationsthereby lowering runtime in line with “green AI” principles. Thus, LM-LRSGA provides a step toward methods that are not only theoretically and numerically robust, but also environmentally sustainable.

Finally, quantifying energy savings from optimization algorithms—both in computational cost and carbon footprint—has become an active research topic~\cite{colao2025optimizer,foglia2024halpernsgd}. In that spirit, future studies will assess the effective carbon footprint reduction enabled by LM-LRSGA, providing a unified view of its theoretical, computational, and environmental impact.

\section*{Acknowledgment(s)}
This publication was partially funded by the PhD program in Mathematics and Computer Science at University of Calabria, Cycle XXXVIII with the support of a scholarship financed by DM 351/2022 (CUP H23C22000440007), based on the NRPP funded by the European Union.

\bibliographystyle{unsrtnat}
\bibliography{references}  %%% Uncomment this line and comment out the ``thebibliography'' section below to use the external .bib file (using bibtex) .





% \end{thebibliography}

\end{document}